\documentclass[11pt]{amsart}

\theoremstyle{plain}
\newtheorem{theorem}                {Theorem}      [section]
\newtheorem{proposition}  [theorem]  {Proposition}
\newtheorem{corollary}    [theorem]  {Corollary}

\theoremstyle{definition}

\newtheorem{remark}       [theorem]  {Remark}

\DeclareMathOperator{\trace}{trace} 
\DeclareMathOperator{\grad}{grad}

\DeclareMathOperator{\vol}{Vol}

\DeclareMathOperator{\Span}{span}
 
\DeclareMathOperator{\cst}{constant}

\numberwithin{equation}{section}

\begin{document}

\title[Biharmonic integral $\mathcal{C}$-parallel submanifolds]
{Biharmonic integral $\mathcal{C}$-parallel submanifolds in
$7$-dimensional Sasakian space forms}

\author{D.~Fetcu}
\author{C.~Oniciuc}

\address{Department of Mathematics\\
"Gh. Asachi" Technical University of Iasi\\
Bd. Carol I no. 11 \\
700506 Iasi, Romania} \email{dfetcu@math.tuiasi.ro}

\curraddr{IMPA\\ Estrada Dona Castorina\\ 110, 22460-320 Rio de
Janeiro, Brasil} \email{dorel@impa.br}

\address{Faculty of Mathematics\\ "Al.I. Cuza" University of Iasi\\
Bd. Carol I no. 11 \\
700506 Iasi, Romania} \email{oniciucc@uaic.ro}

\thanks{The first author was supported by a Post-Doctoral Fellowship "P\'os-Doutorado J\'unior
(PDJ)" offered by CNPq Brazil.}

\begin{abstract} We find the characterization of maximum dimensional
proper-bi\-har\-mo\-nic integral $\mathcal{C}$-parallel submanifolds
of a Sasakian space form and then classify such submanifolds in a
$7$-dimensional Sasakian space form. Working in the sphere
$\mathbb{S}^7$ we explicitly find all $3$-dimensional
proper-biharmonic integral $C$-parallel submanifolds. We also
determine the proper-biharmonic parallel Lagrangian submanifolds of
$\mathbb{C}P^3$.
\end{abstract}

\date{}



\maketitle

\section{Introduction}
\setcounter{equation}{0}

Although, according to its age, the study of biharmonic maps could
be considered a rather old problem, in fact the literature on this
subject experienced an intensive growth in the last decade.

Suggested in 1964, by Eells and Sampson in their famous paper
\cite{JEJS}, as a natural generalization of \textit{harmonic maps}
$\psi:(M,g)\rightarrow(N,h)$ between Riemannian manifolds, which are
critical points of the \textit{energy functional}
$$
E(\psi)=\frac{1}{2}\int_{M}|d\psi|^{2}\ v_{g},
$$
the \textit{biharmonic maps} are critical points of the
\textit{bienergy functional}
$$
E_{2}(\psi)=\frac{1}{2}\int_{M}|\tau(\psi)|^{2}\ v_{g},
$$
where $\tau(\psi)=\trace\nabla d\psi$ is the tension field that
vanishes for harmonic maps. The Euler-Lagrange equation for the
bienergy functional was derived by Jiang in 1986 (see \cite{GYJ}):
$$
\begin{array}{cl}
\tau_{2}(\psi)&=-\Delta\tau(\psi)-\trace\ R^{N}(d\psi,\tau(\psi))d\psi\\
&=0\end{array}
$$
where $\tau_{2}(\psi)$ is the \textit{bitension field} of $\psi$.
Since any harmonic map is biharmonic, we are interested in
non-harmonic biharmonic maps, which are called
\textit{proper-biharmonic}.

An important case of biharmonic maps is represented by the
biharmonic Riemannian immersions, or biharmonic submanifolds, i.e.
submanifolds for which the inclusion map is biharmonic. In Euclidean
spaces the biharmonic submanifolds are the same as those defined by
Chen in \cite{BYC}, as they are characterized by the equation
$\Delta H=0$, where $H$ is the mean curvature vector field and
$\Delta$ is the rough Laplacian.

Pursuing the founding of proper-biharmonic submanifolds in
Riemannian manifolds the attention was first focused on space forms,
and classification results in this context were obtained, for
example, in \cite{ABSMCO, RCSMCO, BYC, ID}. More recently such
results were also found in spaces of non-constant sectional
curvature (see, for example, \cite{RCSMPP, TIJIHU, YLO1, YLO2, WZ}).

A different and active research direction is the study of
proper-biharmonic submanifolds in pseudo-Riemannian manifolds (see,
for example, \cite{AAFDGKVP, AAFDGK, BYC2}).

During the efforts of studying the biharmonic submanifolds in space
forms, the Euclidean spheres proved to be a very giving environment
for obtaining examples and classification results (see \cite{AB} for
detailed proofs). Then, the fact that odd-dimensional spheres can be
thought as a class of Sasakian space forms (which do not have
constant sectional curvature, in general) led to the idea that the
next step would be the study of biharmonic submanifolds in Sasakian
space forms. Following this direction, in \cite{JI} were classified
the proper-biharmonic Legendre curves and Hopf cylinders in a
$3$-dimensional Sasakian space form, whilst in \cite{DFCO3} their
parametric equations were found. In \cite{DFCO} all
proper-biharmonic Legendre curves in any dimensional Sasakian space
forms were classified, and it was provided a method to obtain
proper-biharmonic anti-invariant submanifolds from proper-biharmonic
integral submanifolds. Also, classification results for
proper-biharmonic hypersurfaces were obtained in \cite{DFCO4}.

The goals of our paper are to characterize the maximum dimensional
proper-biharmonic integral, and integral $\mathcal{C}$-parallel,
submanifolds in a Sasakian space form, and then to use these results
in order to obtain the $3$-dimensional proper-biharmonic integral
$\mathcal{C}$-parallel submanifolds of a $7$-dimensional Sasakian
space form. The paper is organized as follows. In Section $2$ we
briefly recall some general facts on Sasakian space forms with a
special emphasis on the notion of integral $\mathcal{C}$-parallel
submanifolds, and also present some old and new results concerning
the proper-biharmonic submanifolds in odd-dimensional spheres.
Section $3$ is devoted to the study of the biharmonicity of maximum
dimensional integral submanifolds in a Sasakian space form. We
obtain the necessary and sufficient conditions for such a
submanifold to be biharmonic, prove some non-existence results and
find the characterization of proper-biharmonic integral
$\mathcal{C}$-parallel submanifolds of maximum dimension. In Section
$4$ we classify all $3$-dimensional proper-biharmonic integral
$\mathcal{C}$-parallel submanifolds in a $7$-dimensional Sasakian
space form, whilst in Section $5$ we find these submanifolds in the
$7$-sphere endowed with its canonical and deformed Sasakian
structures introduced by Tanno in \cite{ST}. In the last section we
classify the proper-biharmonic parallel Lagrangian submanifolds of
$\mathbb{C}P^3$ by determining their horizontal lifts, with respect
to the Hopf fibration, in $\mathbb{S}^7(1)$.

For a general account of biharmonic maps see \cite{SMCO} and
\textit{The Bibliography of Biharmonic Maps} \cite{bibl}.

\noindent \textbf{Conventions.} We work in the $C^{\infty}$
category, that means manifolds, metrics, connections and maps are
smooth. The Lie algebra of the vector fields on $M$ is denoted by
$C(TM)$. The manifold $M$ is always assumed to be connected.

\noindent \textbf{Acknowledgements.} The authors wish to thank
Professor David Blair for useful comments and constant
encouragement, and Professor Harold Rosenberg for helpful
discussions.

\section{Preliminaries}
\setcounter{equation}{0}

\subsection{Integral $\mathcal{C}$-parallel submanifolds of a Sasakian manifold}
A \textit{contact metric structure} on an odd-dimensional manifold
$N^{2n+1}$ is given by $(\varphi,\xi,\eta,g)$, where $\varphi$ is
a tensor field of type $(1,1)$ on $N$, $\xi$ is a vector field,
$\eta$ is a 1-form and $g$ is a Riemannian metric such that
$$
\varphi^{2}=-I+\eta\otimes\xi,\quad \eta(\xi)=1
$$
and
$$
g(\varphi U,\varphi V)=g(U,V)-\eta(U)\eta(V),\quad g(U,\varphi
V)=d\eta(U,V),\quad \forall U,V\in C(TN).
$$

\noindent A contact metric structure $(\varphi,\xi,\eta,g)$ is
called {\it normal} if
$$
N_{\varphi}+2d\eta\otimes\xi=0,
$$
where
$$
N_{\varphi}(U,V)=[\varphi U,\varphi V]-\varphi \lbrack \varphi
U,V]-\varphi \lbrack U,\varphi V]+\varphi ^{2}[U,V],\quad \forall
U,V\in C(TN),
$$
is the Nijenhuis tensor field of $\varphi$.

A contact metric manifold $(N,\varphi,\xi,\eta,g)$ is
\textit{regular} if for any point $p\in N$ there exists a cubic
neighborhood such that any integral curve of $\xi$ passes through it
at most once; and it is \textit{strictly regular} if all integral
curves of $\xi$ are homeomorphic to each other.

A contact metric manifold $(N,\varphi,\xi,\eta,g)$ is a
\textit{Sasakian manifold} if it is normal or, equivalently, if
$$
(\nabla^N_{U}\varphi)(V)=g(U,V)\xi-\eta(V)U,\quad\forall U,V\in
C(TN),
$$
where $\nabla^N$ is the Levi-Civita connection on $(N,g)$. We shall
often use in our paper the formula $\nabla^N_U\xi=-\varphi U$, which
holds on a Sasakian manifold.

\noindent Let $(N,\varphi,\xi,\eta,g)$ be a Sasakian manifold. The
sectional curvature of a 2-plane generated by $U$ and $\varphi U$,
where $U$ is a unit vector orthogonal to $\xi$, is called
\textit{$\varphi$-sectional curvature} determined by $U$. A Sasakian
manifold with constant $\varphi$-sectional curvature $c$ is called a
\textit{Sasakian space form} and is denoted by $N(c)$. The curvature
tensor field of a Sasakian space form $N(c)$ is given by
$$
\begin{array}{ll}
R^N(U,V)W=&\frac{c+3}{4}\{g(W,V)U-g(W,U)V\}+\frac{c-1}{4}\{\eta(W)\eta(U)V\\ \\
&-\eta(W)\eta(V)U+g(W,U)\eta(V)\xi-g(W,V)\eta(U)\xi\\
\\&+g(W,\varphi V)\varphi U-g(W,\varphi U)\varphi
V+2g(U,\varphi V)\varphi W\}.
\end{array}
$$

\noindent The classification of the complete, simply connected
Sasakian space forms $N(c)$ was given in \cite{ST}. Thus, if $c=1$
then $N(1)$ is isometric to the unit sphere $\mathbb{S}^{2n+1}$
endowed with its canonical Sasakian structure and if $c>-3$ then
$N(c)$ is isometric to $\mathbb{S}^{2n+1}$ endowed with the deformed
Sasakian structure introduced by Tanno in \cite{ST}, which we
present below.

Let $\mathbb{S}^{2n+1}=\{z\in\mathbb{C}^{n+1}: \vert z\vert=1\}$ be
the unit $(2n+1)$-dimensional Euclidean sphere. Consider the
following structure tensor fields on $\mathbb{S}^{2n+1}$:
$\xi_{0}=-\mathcal{J}z$, for each $z\in \mathbb{S}^{2n+1}$, where
$\mathcal{J}$ is the usual complex structure on $\mathbb{C}^{n+1}$
defined by
$$
\mathcal{J}z=(-y^{1},...,-y^{n+1},x^{1},...,x^{n+1}),
$$
for $z=(x^{1},...,x^{n+1},y^{1},...,y^{n+1})$, and
$\varphi_{0}=s\circ \mathcal{J}$, where $s:T_{z}\mathbb{C}^{n+1}\to
T_{z}\mathbb{S}^{2n+1}$ denotes the orthogonal projection. Equipped
with these tensors and the standard metric $g_0$, the sphere
$\mathbb{S}^{2n+1}$ becomes a Sasakian space form with
$\varphi_{0}$-sectional curvature equal to $1$, denoted by
$\mathbb{S}^{2n+1}(1)$.

\noindent Now, consider the deformed Sasakian structure on
$\mathbb{S}^{2n+1}$,
$$
\eta=a\eta_{0},\quad\xi=\frac{1}{a}\xi_{0},\quad
\varphi=\varphi_{0},\quad g=a g_{0}+a(a-1)\eta_{0}\otimes\eta_{0},
$$
where $a$ is a positive constant. The structure
$(\varphi,\xi,\eta,g)$ is still a Sasakian structure and
$(\mathbb{S}^{2n+1},\varphi,\xi,\eta,g)$ is a Sasakian space form
with constant $\varphi$-sectional curvature $c=\frac{4}{a}-3>-3$,
denoted by $\mathbb{S}^{2n+1}(c)$ (see also \cite{DB}).

A submanifold $M^m$ of a Sasakian manifold
$(N^{2n+1},\varphi,\xi,\eta,g)$ is called an {\it integral
submanifold} if $\eta(X)=0$ for any vector field $X$ tangent to
$M$. We have $\varphi(TM)\subset NM$ and $m\leq n$, where $TM$ and
$NM$ are the tangent bundle and the normal bundle of $M$,
respectively. Moreover, for $m=n$, one gets $\varphi(NM)=TM$. If
we denote by $B$ the second fundamental form of $M$ then, by a
straightforward computation, one obtains the following relation
$$
g(\varphi Z,B(X,Y))=g(\varphi Y,B(X,Z)),
$$
for any vector fields $X,Y$ and $Z$ tangent to $M$ (see also
\cite{CBDBTK}). We also note that $A_\xi=0$, where $A$ is the shape
operator of $M$ (see \cite{DB}).

A submanifold $\widetilde M$ of $N$ is said to be
\textit{anti-invariant} if it is tangent to $\xi$ and
$\varphi(T\widetilde M)\subset N\widetilde M$.

Next, we shall recall the notion of an integral
$\mathcal{C}$-parallel submanifold of a Sasakian manifold (see, for
example, \cite{CBDBTK}). Let $M^m$ be an integral submanifold of a
Sasakian manifold $(N^{2n+1},\varphi,\xi,\eta,g)$. Then $M$ is said
to be \textit{integral $\mathcal{C}$-parallel} if $\nabla^{\perp}B$
is parallel to the characteristic vector field $\xi$, where $B$ is
the second fundamental form of $M$ and $\nabla^{\perp}B$ is given by
$$
(\nabla^{\perp}
B)(X,Y,Z)=\nabla^{\perp}_{X}B(Y,Z)-B(\nabla_{X}Y,Z)-B(Y,\nabla_{X}Z)
$$
for any vector fields $X,Y,Z$ tangent to $M$, $\nabla^{\perp}$ and
$\nabla$ being the normal connection and the Levi-Civita connection
on $M$, respectively. This means $(\nabla^{\perp}B)(X,Y,Z)=g(\varphi
X,B(Y,Z))\xi$. If we denote $S(X,Y,Z)=g(\varphi X,B(Y,Z))$ then $S$
is a totally symmetric tensor field of type $(0,3)$ on $M$.

In general, a submanifold $M$ of $N$ is called \textit{parallel} if
$\nabla^{\perp}B=0$.

The following two results shall be used latter in this paper and,
for the sake of completeness, we also provide their proofs.

\begin{proposition}\label{prop00} If the mean curvature vector field $H$ of an
integral submanifold $M^{n}$ of a Sasakian manifold
$(N^{2n+1},\varphi,\xi,\eta,g)$ is parallel then $M^{n}$ is
minimal.
\end{proposition}
\begin{proof}
Let $X,Y$ be two vector fields tangent to $M$. Since
$$
g(B(X,Y),\xi)=g(\nabla^{N}_{X}Y,\xi)=-g(Y,\nabla^N_X\xi)=g(Y,\varphi
X)=0
$$
we have $B(X,Y)\in \varphi(TM)$ and, in particular,
$H\in\varphi(TM)$. Then
$$
g(\nabla^{\perp}_{X}H,\xi)=g(\nabla^{N}_{X}H,\xi)=-g(H,\nabla^N_X\xi)=g(H,\varphi
X).
$$
Thus, if $\nabla^{\perp}H=0$ it follows that $g(H,\varphi X)=0$
for any vector field $X$ tangent to $M$, and this means $H=0$.
\end{proof}

\begin{proposition}\label{prop0} Let $(N^{2n+1},\varphi,\xi,\eta,g)$ be
a Sasakian manifold and $M^n$ be an integral $\mathcal{C}$-parallel
submanifold with mean curvature vector field $H$. The following
hold:
\begin{enumerate}
\item $\nabla^{\perp}_{X}H=g(H,\varphi X)\xi$, for any vector
field $X$ tangent to $M$, i.e. $H$ is $\mathcal{C}$-parallel;

\item $\Delta^{\perp}H=H$;

\item the mean curvature $|H|$ is constant.
\end{enumerate}
\end{proposition}
\begin{proof} Consider
$\{X_{i}\}_{i=1}^{n}$ to be a local geodesic frame at $p\in M$.
Then we have at $p$
$$
(\nabla^{\perp}
B)(X_{i},X_{j},X_{j})=\nabla^{\perp}_{X_{i}}B(X_{j},X_{j})=g(B(X_{j},X_{j}),\varphi
X_{i})\xi
$$
and, by summing after $j=\overline{1,n}$, we obtain
$\nabla^{\perp}_{X_{i}}H=g(H,\varphi X_{i})\xi$.

\noindent Next, as $\nabla^N_X\xi=-\varphi X$, from the Weingarten
equation we get $A_\xi=0$, where $A_{\xi}$ is the shape operator of
$M$ corresponding to $\xi$, and
$\nabla^{\perp}_X\xi=\nabla^N_X\xi=-\varphi X$. Thus
$$
\begin{array}{ll}
\Delta^{\perp}H&=-\sum_{i=1}^{n}\nabla_{X_{i}}^{\perp}\nabla_{X_{i}}^{\perp}H=
-\sum_{i=1}^{n}\nabla_{X_{i}}^{\perp}(g(H,\varphi X_{i})\xi)\\
\\&= -\sum_{i=1}^{n}X_{i}(g(H,\varphi
X_{i}))\xi-\sum_{i=1}^{n}(g(H,\varphi
X_{i}))\nabla_{X_{i}}^{N}\xi\\
\\&=-\sum_{i=1}^{n}X_{i}(g(H,\varphi
X_{i}))\xi+\sum_{i=1}^{n}(g(H,\varphi X_{i}))\varphi X_{i}\\
\\&=-\sum_{i=1}^{n}X_{i}(g(H,\varphi X_{i}))\xi+H.
\end{array}
$$
But, since $\nabla_{X_{i}}^{N}\varphi
X_{i}=\varphi\nabla^{N}_{X_{i}}X_{i}+\xi$, it results
$$
\begin{array}{ll}
X_{i}(g(H,\varphi X_{i}))&=g(\nabla_{X_{i}}^{N}H,\varphi
X_{i})+g(H,\varphi\nabla^{N}_{X_{i}}X_{i}+\xi)\\
\\&=g(-A_{H}X_{i}+\nabla^{\perp}_{X_{i}}H,\varphi X_{i})+g(H,\varphi
B(X_{i},X_{i}))\\ \\&=0.
\end{array}
$$

\noindent We have just proved that $\Delta^{\perp}H=H$.

\noindent Finally, we have
$$
X(|H|^{2})=2g(H,\nabla^{\perp}_{X}H)=2g(H,\varphi X)g(H,\xi)=0
$$
for any vector field $X$ tangent to $M$.
Consequently, it follows $|H|=\cst$.
\end{proof}

\subsection{Biharmonic submanifolds in $\mathbb{S}^{2n+1}(1)$}

We shall recall first the notion of Frenet curve of osculating order
$r$ as it is presented, for example, in \cite{HN}. Let $(M^{m},g)$
be a Riemannian manifold and $\Gamma:I\to M$ a curve parametrized by
arc length, that is $|\Gamma'|=1$. Then $\Gamma$ is called a
\textit{Frenet curve of osculating order $r$}, $1\leq r\leq m$, if
for all $s\in I$ its higher order derivatives
$$
\Gamma'(s)=(\nabla^0_{\Gamma'}\Gamma')(s),\quad(\nabla_{\Gamma'}\Gamma')(s),\quad
\dots,\quad(\nabla^{r-1}_{\Gamma'}\Gamma')(s)
$$
are linearly independent but
$$
\Gamma'(s)=(\nabla^0_{\Gamma'}\Gamma')(s),\quad(\nabla_{\Gamma'}\Gamma')(s),\quad
\dots,\quad(\nabla^{r-1}_{\Gamma'}\Gamma')(s),\quad(\nabla^r_{\Gamma'}\Gamma')(s)
$$
are linearly dependent in $T_{\Gamma(s)}M$. Then there exist
unique orthonormal vector fields $E_{1},E_{2},...,E_{r}$ along
$\Gamma$ such that
$$
\nabla_{T}E_{1}=\kappa_{1}E_{2},\quad
\nabla_{T}E_{2}=-\kappa_{1}E_{1}+\kappa_{2}E_{3},...,\nabla_{T}E_{r}=-\kappa_{r-1}E_{r-1}
$$
where $E_{1}=\Gamma'=T$ and $\kappa_{1},...,\kappa_{r-1}$ are
positive functions on $I$.

\begin{remark}A geodesic is a Frenet curve of osculating order $1$;
a \textit{circle} is a Frenet curve of osculating order $2$ with
$\kappa_{1}=\cst$; a \textit{helix of order r}, $r\geq 3$, is a
Frenet curve of osculating order $r$ with
$\kappa_{1},...,\kappa_{r-1}$ constants; a helix of order $3$ is
simply called a \textit{helix}.
\end{remark}

In \cite{JI} Inoguchi proved that there are no proper-biharmonic
Legendre curves in $\mathbb{S}^3(1)$ whilst in \cite{DFCO} we found
the parametric equations of all proper-biharmonic Legendre curves in
$\mathbb{S}^{2n+1}(1)$, $n\geq 2$. These curves are given by the
following

\begin{theorem} [\cite{DFCO}]\label{curv2s2n+1,1}
Let $\Gamma:I\to
(\mathbb{S}^{2n+1},\varphi_{0},\xi_{0},\eta_{0},g_{0})$, $n\geq 2$,
be a proper-biharmonic Legendre curve parametrized by arc length.
Then the parametric equation of $\Gamma$ in the Euclidean space
$(\mathbb{R}^{2n+2},\langle,\rangle)$, is either
$$
\Gamma(s)=\frac{1}{\sqrt{2}}\cos(\sqrt{2}s)e_{1}+
\frac{1}{\sqrt{2}}\sin(\sqrt{2}s)e_{2}+\frac{1}{\sqrt{2}}e_{3}
$$
where $\{e_{i},\mathcal{J}e_{j}\}_{i,j=1}^{3}$ are constant unit
vectors orthogonal to one another, or
$$
\Gamma(s)=\frac{1}{\sqrt{2}}\cos(As)e_{1}+\frac{1}{\sqrt{2}}\sin(As)e_{2}+
\frac{1}{\sqrt{2}}\cos(Bs)e_{3}+\frac{1}{\sqrt{2}}\sin(Bs)e_{4},
$$
where
$$
A=\sqrt{1+\kappa_{1}},\quad B=\sqrt{1-\kappa_{1}},\quad
\kappa_{1}\in(0,1)
$$
and $\{e_{i}\}_{i=1}^{4}$ are constant unit vectors orthogonal to
one another, satisfying
$$
\langle e_{1},\mathcal{J}e_{3}\rangle=\langle
e_{1},\mathcal{J}e_{4}\rangle=\langle
e_{2},\mathcal{J}e_{3}\rangle=\langle
e_{2},\mathcal{J}e_{4}\rangle=0,\quad A\langle
e_{1},\mathcal{J}e_{2}\rangle+B\langle
e_{3},\mathcal{J}e_{4}\rangle=0.
$$
\end{theorem}

\begin{remark}\label{remark1}
We note that if $\Gamma$ is a proper-biharmonic Legendre circle,
then $E_2\perp\varphi T$ and $n\geq 3$. If $\Gamma$ is a
proper-biharmonic Legendre helix, then $g_0(E_2,\varphi T)=-A\langle
e_1, \mathcal{J}e_2\rangle$ and we have two cases: either
$E_2\perp\varphi T$ and then
$\{e_{i},\mathcal{J}e_{j}\}_{i,j=1}^{4}$ is an orthonormal system in
$\mathbb{R}^{2n+2}$, so $n\geq 3$, or $g_0(E_2,\varphi T)\neq 0$
and, in this case, $g_0(E_2,\varphi T)\in(-1,1)\setminus\{0\}$. We
also observe that $\varphi T$ cannot be parallel to $E_2$. When
$g_0(E_2,\varphi T)\neq 0$ and $n\geq 3$ the first four vectors (for
example) in the canonical basis of the Euclidean space
$\mathbb{R}^{2n+2}$ satisfy the conditions of Theorem
\ref{curv2s2n+1,1}, whilst for $n=2$ we can obtain four vectors
$\{e_1,e_2,e_3,e_4\}$ satisfying these conditions in the following
way. We consider constant unit vectors $e_1$, $e_3$ and $f$ in
$\mathbb{R}^6$ such that
$\{e_1,e_3,f,\mathcal{J}e_1,\mathcal{J}e_3,\mathcal{J}f\}$ is a
$\mathcal{J}$-basis. Then, by a straightforward computation, it
follows that the vectors $e_2$ and $e_4$ have to be given by
$$
e_2=\mp\frac{B}{A}\mathcal{J}e_1+\alpha_1f+\alpha_2\mathcal{J}f,\quad
e_4=\pm\mathcal{J}e_3,
$$
where $\alpha_1$ and $\alpha_2$ are constants such that
$\alpha_1^2+\alpha_2^2=1-\frac{B^2}{A^2}=\frac{2\kappa_1}{A^2}$. As
a concrete example, we can start with the following vectors in
$\mathbb{R}^6$:
$$
e_1=(1,0,0,0,0,0),\quad e_3=(0,0,1,0,0,0),\quad f=(0,1,0,0,0,0)
$$
and obtain
$$
e_2=\Big(0,\alpha_1,0,-\frac{B}{A},\alpha_2,0\Big),\quad
e_4=(0,0,0,0,0,1),
$$
where $\alpha_1^2+\alpha_2^2=1-\frac{B^2}{A^2}$.
\end{remark}

The classification of all proper-biharmonic Legendre curves in a
Sasakian space form $N^{2n+1}(c)$ was given in \cite{DFCO}. This
classification is invariant under an isometry $\Psi$ of $N$ which
preserves $\xi$ (or, equivalently, $\Psi$ is $\varphi$-holomorphic).

In order to find higher dimensional proper-biharmonic submanifolds
in a Sasakian space form we gave the following

\begin{theorem}[\cite{DFCO}] \label{teorema1}
Let $(N^{2n+1},\varphi,\xi,\eta,g)$ be a strictly regular Sasakian
space form with constant $\varphi$-sectional curvature $c$ and let
${\bf i}:M\to N$ be an $m$-dimensional integral submanifold of $N$,
$1\leq m\leq n$. Consider the cylinder
$$
F:\widetilde{M}=I\times M\to N,\quad
F(t,p)=\phi_{t}(p)=\phi_{p}(t),
$$
where $I=\mathbb{S}^{1}$ or $I=\mathbb{R}$ and $\{\phi_{t}\}_{t\in
I}$ is the flow of the vector field $\xi$. Then
$F:(\widetilde{M},\widetilde{g}=dt^{2}+{\bf i}^{\ast}g)\to N$ is an
anti-invariant Riemannian immersion, and is proper-biharmonic if and
only if $M$ is a proper-biharmonic submanifold of $N$.
\end{theorem}

Working with anti-invariant submanifolds rather cylinders, we can
state the following (known) result

\begin{proposition} Let $\widetilde{M}^{m+1}$ be an anti-invariant submanifold
of the strictly regular Sasakian space form $N^{2n+1}(c)$, $1\leq
m\leq n$, invariant under the flow-action of the characteristic
vector field $\xi$. Then $\widetilde{M}$ is locally isometric to
$I\times M^m$, where $M^m$ is an integral submanifold of $N$.
Moreover, we have
\begin{enumerate}
\item $\widetilde{M}$ is proper-biharmonic if and only if $M$ is
proper-biharmonic in $N$;
\item if $m=n$, then $\widetilde{M}$ is parallel if and only if
$M$ is $\mathcal{C}$-parallel;
\item if $m=n$, then $\widetilde{M}$ has parallel mean curvature
vector field if and only if $M$ satisfies
$\nabla^{\perp}H\parallel\xi$.
\end{enumerate}
\end{proposition}

\begin{proof} The restriction $\xi_{/\widetilde{M}}$ of the characteristic vector field
$\xi$ to $\widetilde{M}$ is a Killing tangent vector field on
$\widetilde{M}$. Since $\widetilde{M}$ is anti-invariant, the
horizontal distribution defined on $\widetilde{M}$ is integrable.
Let $p\in\widetilde{M}$ be an arbitrary point and $M$ a small enough
integral submanifold of the horizontal distribution on
$\widetilde{M}$ such that $p\in M$. Then $F:I\times M\rightarrow
F(I\times M)\subset\widetilde M$, $F(t,p)=\phi_t(p)$, is an
isometry. As $M$ is an integral submanifold of the horizontal
distribution on $\widetilde{M}$, it is an integral submanifold of
$N$.

The item $(1)$ follows immediately from Theorem \ref{teorema1}, and
$(2)$ and $(3)$ are known and can be checked by straightforward
computations.
\end{proof}

As a surface in a strictly regular Sasakian space form which is
invariant under the flow-action of the characteristic vector field
is also anti-invariant, we have

\begin{corollary}\label{teorema2} Let $\widetilde{M}^2$ be a
surface of $N^{2n+1}(c)$ invariant under the flow-action of the
characteristic vector field $\xi$. Then $\widetilde{M}$ is locally
isometric to $I\times\Gamma$, where $\Gamma$ is a Legendre curve in
$N$ and, moreover, it is proper-biharmonic if and only if $\Gamma$
is proper-biharmonic in $N$.
\end{corollary}

Now, consider $\widetilde{M}^2$ a surface of $N^{2n+1}(c)$
invariant under the flow-action of the characteristic vector field
$\xi$ and let $T=\Gamma'$ and $E_2$ be the first two vector fields
defined by the Frenet equations of the above Legendre curve
$\Gamma$. As $\nabla^{F}_{\partial/\partial
t}\tau(F)=-\varphi(\tau(F))$, where $\nabla^{F}$ is the pull-back
connection determined by the Levi-Civita connection on $N$, we can
prove

\begin{proposition}\label{propsurfinv} Let $\widetilde{M}^2$ be a
proper-biharmonic surface of $N^{2n+1}(c)$ invariant under the
flow-action of the characteristic vector field $\xi$. Then
$\widetilde{M}$ has parallel mean curvature vector field if and only
if $c>1$ and $\varphi T\parallel E_2$.
\end{proposition}

From Proposition \ref{propsurfinv} it results

\begin{corollary} The proper-biharmonic surfaces of $\mathbb{S}^{2n+1}(1)$ invariant under
the flow-action of the characteristic vector field $\xi_0$ are not
of parallel mean curvature vector field.
\end{corollary}

We shall see that we do have examples of maximum dimensional
proper-bi\-har\-mo\-nic anti-invariant submanifolds of
$\mathbb{S}^{2n+1}(1)$, invariant under the flow-action of
$\xi_0$, which have parallel mean curvature vector field.

In \cite{TS} the parametric equations of all proper-biharmonic
integral surfaces in $\mathbb{S}^{5}(1)$ were obtained. Up to an
isometry of $\mathbb{S}^5(1)$ which preserves  $\xi_{0}$, we have
only one proper-biharmonic integral surface given by
$$
x(u,v)=\frac{1}{\sqrt{2}}(\exp(\mathrm{i}u),\mathrm{i}\exp(-\mathrm{i}u)
\sin(\sqrt{2}v),\mathrm{i}\exp(-\mathrm{i}u)\cos(\sqrt{2}v)).
$$
The map $x$ induces a proper-biharmonic Riemannian embedding from
the $2$-di\-men\-si\-o\-nal torus
$\mathcal{T}^2=\mathbb{R}^2/\Lambda$ into $\mathbb{S}^5$, where
$\Lambda$ is the lattice generated by the vectors $(2\pi,0)$ and
$(0,\sqrt{2}\pi)$.

\begin{remark} We recall that an isometric immersion
$x:M\rightarrow\mathbb{R}^{n+1}$ of a compact manifold is said to be
of \textit{$k$-type} if its spectral decomposition contains exactly
$k$ non-zero terms excepting the center of mass
$x_0=\frac{1}{\vol(M)}\int_Mx\ v_g$. When $x_0=0$, the submanifold
is called \textit{mass-symmetric} (see ~\cite{BYC}). It was proved
in \cite{ABSMCO,ABSMCO2} that, in general, a proper-biharmonic
compact constant mean curvature submanifold $M^m$ of $\mathbb{S}^n$
is either a $1$-type submanifold of $\mathbb{R}^{n+1}$ with center
of mass of norm equal to $\frac{1}{\sqrt{2}}$, or is a
mass-symmetric $2$-type submanifold of $\mathbb{R}^{n+1}$. Now,
using Theorem $3.5$ in \cite{CBDB}, where all mass-symmetric
$2$-type integral surfaces in $\mathbb{S}^5(1)$ were determined, and
Proposition $4.1$ in \cite{RCSMCO}, the result in \cite{TS} can be
(partially) reobtained.
\end{remark}

Further, we consider the cylinder over $x$ and we recover the result
in \cite{KARECMTS}: up to an isometry which preserves $\xi_0$, we
have only one $3$-dimensional proper-biharmonic anti-invariant
submanifold of $\mathbb{S}^5(1)$ invariant under the flow-action of
$\xi_0$,
$$
y(t,u,v)=\exp(-\mathrm{i}t)x(u,v).
$$
The map $y$ is a proper-biharmonic Riemannian immersion with
parallel mean curvature vector field and induces a proper-biharmonic
Riemannian immersion from the $3$-dimensional torus
$\mathcal{T}^3=\mathbb{R}^3/\Lambda$ into $\mathbb{S}^5$, where
$\Lambda$ is the lattice generated by the vectors $(2\pi,0,0)$,
$(0,2\pi,0)$ and $(0,0,\sqrt{2}\pi)$. Moreover, a closer look shows
that $y$ factorizes to a proper-biharmonic Riemannian embedding in
$\mathbb{S}^5$ and its image is the Riemannian product between three
Euclidean circles, one of radius $\frac{1}{\sqrt{2}}$ and each of
the other two of radius $\frac{1}{2}$. Indeed, we may consider the
orthogonal transformation of $\mathbb{R}^3$ given by
$$
T(t,u,v)=\Big(\frac{-t+u}{\sqrt{2}},\frac{-t-u}{\sqrt{2}},v\Big)=(t',u',v')
$$
and the map $y$ becomes
$$
y_1(t',u',v')=\frac{1}{\sqrt{2}}
(\exp({\mathrm{i}\sqrt{2}t'}),\mathrm{i}\exp({\mathrm{i}\sqrt{2}u'})\sin(\sqrt{2}v'),
\mathrm{i}\exp({\mathrm{i}\sqrt{2}u'})\cos(\sqrt{2}v')).
$$
Then, acting with an appropriate holomorphic isometry of
$\mathbb{C}^4$, $y_1$ becomes
$$
y_2(t',u',v')=\Big(\frac{1}{\sqrt{2}}
\exp({\mathrm{i}\sqrt{2}t'}),\frac{1}{2}\exp({\mathrm{i}(u'-v'})),
\frac{1}{2}\exp({\mathrm{i}(u'+v'}))\Big)
$$
and, further, an obvious orthogonal transformation of the domain
leads to the desired results.

\section{Biharmonic integral submanifolds of maximum dimension in
Sasakian space forms} \setcounter{equation}{0}

Let $(N^{2n+1},\varphi,\xi,\eta,g)$ be a Sasakian space form with
constant $\varphi$-sectional curvature $c$, and and $M^n$ be an
$n$-dimensional integral submanifold of $N$. We shall denote by
$B$, $A$ and $H$ the second fundamental form of $M$ in $N$, the
shape operator and the mean curvature vector field, respectively.
By $\nabla^{\perp}$ and $\Delta^{\perp}$ we shall denote the
connection and the Laplacian in the normal bundle. We have

\begin{theorem}\label{t0}
The integral submanifold ${\bf i}:M^{n}\to N^{2n+1}$ is biharmonic
if and only if
\begin{equation}\label{ec3.0}
\begin{cases}\Delta^{\perp}H+\trace
B(\cdot,A_{H}\cdot)
-\frac{c(n+3)+3n-3}{4}H=0\\
4\trace
A_{\nabla^{\perp}_{(\cdot)}H}(\cdot)+n\grad(|H|^{2})=0.\end{cases}
\end{equation}
\end{theorem}
\begin{proof} Let us denote by $\nabla^{N}$, $\nabla$ the
Levi-Civita connections on $N$ and $M$, respectively. Consider
$\{X_{i}\}_{i=1}^{n}$ to be a local geodesic frame at $p\in M$.
Then, since $\tau({\bf i})=nH$, we have at $p$
\begin{equation}\label{ec3.2}
\begin{array}{ll}
\tau_{2}({\bf i})&=-\Delta\tau({\bf i})-\trace R^{N}(d{\bf
i},\tau({\bf i}))d{\bf i}\\ \\
&=n\{\sum_{i=1}^{n}\nabla_{X_{i}}^{N}\nabla_{X_{i}}^{N}H-\sum_{i=1}^{n}
R^{N}(X_{i},H)X_{i}\}.\end{array}
\end{equation}

We recall the Weingarten equation, around $p$,
$$
\nabla^{N}_{X_{i}}H=\nabla_{X_{i}}^{\perp}H-A_{H}(X_{i})
$$
and, using the Weingarten and Gauss equations,
$$
\nabla^{N}_{X_{i}}\nabla^{N}_{X_{i}}H=\nabla_{X_{i}}^{\perp}
\nabla_{X_{i}}^{\perp}H-A_{\nabla_{X_{i}}^{\perp}H}(X_{i})
-\nabla_{X_{i}}A_{H}(X_{i})-B(X_{i},A_{H}(X_{i})).
$$
\noindent Thus, at $p$, one obtains
\begin{equation}\label{ec3.3}
\begin{array}{ll}
-\frac{1}{n}\Delta\tau({\bf
i})&=\sum_{i=1}^{n}\nabla^{N}_{X_{i}}\nabla^{N}_{X_{i}}H\\
\\&=-\Delta^{\perp}H-\trace B(\cdot,A_{H}\cdot)-\trace
A_{\nabla^{\perp}_{(\cdot)}H}(\cdot)-\trace \nabla
A_{H}(\cdot,\cdot).
\end{array}
\end{equation}

The next step is to compute $\trace \nabla A_{H}(\cdot,\cdot)$. We
obtain at $p$
$$
\begin{array}{lll}
\trace \nabla
A_{H}(\cdot,\cdot)&=&\sum_{i=1}^{n}\nabla_{X_{i}}A_{H}(X_{i})
=\sum_{i,j=1}^{n}\nabla_{X_{i}}(g(A_{H}(X_{i}),X_{j})X_{j})\\
\\&=&\sum_{i,j=1}^{n}X_{i}(g(A_{H}(X_{i}),X_{j}))X_{j}\\ \\
&=&\sum_{i,j=1}^{n}X_{i}(g(B(X_{j},X_{i}),H))X_{j}\\
\\&=&\sum_{i,j=1}^{n}X_{i}(g(\nabla_{X_{j}}^{N}X_{i},H))X_{j}\\
\\&=&\sum_{i,j=1}^{n}\{g(\nabla_{X_{i}}^{N}\nabla_{X_{j}}^{N}X_{i},H)
+g(\nabla_{X_{j}}^{N}X_{i},\nabla_{X_{i}}^{N}H)\}X_{j}\\
\\&=&\sum_{i,j=1}^{n}g(\nabla_{X_{i}}^{N}\nabla_{X_{j}}^{N}X_{i},H)X_{j}
+\sum_{i,j=1}^{n}g(B(X_{j},X_{i}),\nabla_{X_{i}}^{\perp}H)X_{j}\\
\\&=&\sum_{i,j=1}^{n}g(\nabla_{X_{i}}^{N}\nabla_{X_{j}}^{N}X_{i},H)X_{j}
+\sum_{i,j=1}^{n}g(A_{\nabla_{X_{i}}^{\perp}H}(X_{i}),X_{j})X_{j}\\
\\&=&\sum_{i,j=1}^{n}g(\nabla_{X_{i}}^{N}\nabla_{X_{j}}^{N}X_{i},H)X_{j}+\trace
A_{\nabla_{(\cdot)}^{\perp}H}(\cdot).
\end{array}
$$
Further, using the expression of the curvature tensor field $R^N$,
we have
\begin{equation}\label{ec3.4}
\begin{array}{lll}
\trace \nabla
A_{H}(\cdot,\cdot)&=&\sum_{i,j=1}^{n}g(\nabla_{X_{j}}^{N}\nabla_{X_{i}}^{N}X_{i}
+R^{N}(X_{i},X_{j})X_{i}+\nabla_{[X_{i},X_{j}]}^{N}X_{i},H)X_{j}\\
\\&&+\trace
A_{\nabla_{(\cdot)}^{\perp}H}(\cdot)\\
\\&=&\sum_{i,j=1}^{n}g(\nabla_{X_{j}}^{N}\nabla_{X_{i}}^{N}X_{i},H)X_{j}
+\sum_{i,j=1}^{n}g(R^{N}(X_{i},X_{j})X_{i},H)X_{j}\\ \\&&+\trace
A_{\nabla_{(\cdot)}^{\perp}H}(\cdot).
\end{array}
\end{equation}
\noindent But
\begin{equation}\label{no1}
\begin{array}{lll}
\sum_{i,j=1}^{n}g(\nabla_{X_{j}}^{N}\nabla_{X_{i}}^{N}X_{i},H)X_{j}
&=&\sum_{i,j=1}^{n}g(\nabla_{X_{j}}^{N}B(X_{i},X_{i}),H)X_{j}\\ \\
&&+\sum_{i,j=1}^{n}g(\nabla_{X_{j}}^{N}\nabla_{X_{i}}X_{i},H)X_{j}\\
\\&=&n\sum_{j=1}^{n}g(\nabla^{N}_{X_{j}}H,H)X_{j}\\
\\&&+\sum_{i,j=1}^{n}g(\nabla_{X_{j}}\nabla_{X_{i}}X_{i}
+B(X_{j},\nabla_{X_{i}}X_{i}),H)X_{j}\\
\\&=&\frac{n}{2}\grad(|H|^{2})
\end{array}
\end{equation}
and
\begin{equation}\label{no2}
\begin{array}{lll}
\sum_{i,j=1}^{n}g(R^{N}(X_{i},X_{j})X_{i},H)X_{j}
&=&\sum_{i,j=1}^{n}g(R^{N}(X_{i},H)X_{i},X_{j})X_{j}\\
\\&=&(\trace R^{N}(d{\bf i},H)d{\bf i})^{\top}.\end{array}
\end{equation}

\noindent Replacing (\ref{no1}) and (\ref{no2}) into
(\ref{ec3.4}), we have
$$
\trace\nabla A_{H}(\cdot,\cdot)=\frac{n}{2}\grad(|H|^{2})+(\trace
R^{N}(d{\bf i},H)d{\bf i})^{\top}+\trace
A_{\nabla_{(\cdot)}^{\perp}H}(\cdot)
$$
and therefore
\begin{equation}\label{ec3.5}
\begin{array}{lll}
\trace A_{\nabla_{(\cdot)}^{\perp}H}(\cdot)+\trace\nabla
A_{H}(\cdot,\cdot)&=&2\trace
A_{\nabla_{(\cdot)}^{\perp}H}(\cdot)+\frac{n}{2}\grad(|H|^{2})\\
\\&&+(\trace R^{N}(d{\bf i},H)d{\bf i})^{\top}.\end{array}
\end{equation}

Now, let $\{X_i\}_{i=1}^{n}$ be a local orthonormal frame on $M$.
Then $\{X_{i},\varphi X_{j},\xi\}_{i,j=1}^n$ is a local orthonormal
frame on $N$. By using the expression of the curvature tensor field
and $H\in\Span\{\varphi X_{i}:i=\overline{1,n}\}$ one obtains, after
a straightforward computation,
$$
R^{N}(X_{i},H)X_{i}=-\frac{c+3}{4}H+\frac{3(c-1)}{4}g(\varphi
H,X_{i})\varphi X_{i}.
$$

\noindent Hence
\begin{equation}\label{ec3.6}
\begin{array}{lll} \trace R^{N}(d{\bf i},H)d{\bf
i}&=&\sum_{i=1}^n R^N(X_i,H)X_i\\
\\&=&-\frac{(c+3)n}{4}H+\sum_{i=1}^{n}\frac{3(c-1)}{4}g(\varphi
H,X_{i})\varphi X_{i}\\
\\&=&-\frac{(c+3)n}{4}H-\frac{3(c-1)}{4}H\\
\\&=&-\frac{c(n+3)+3n-3}{4}H,
\end{array}
\end{equation}
which implies $(\trace R^{N}(d{\bf i},H)d{\bf i})^{\top}=0$.

From \eqref{ec3.2}, \eqref{ec3.3}, \eqref{ec3.5} and \eqref{ec3.6}
we have
$$
\begin{array}{lll}
\frac{1}{n}\tau_{2}({\bf i})&=&-\Delta^{\perp}H-\trace
B(\cdot,A_{H}\cdot)+\frac{c(n+3)+3n-3}{4}H\\
\\&&-2\trace
A_{\nabla_{(\cdot)}^{\perp}H}(\cdot)-\frac{n}{2}\grad(|H|^{2}),
\end{array}
$$
and we come to the conclusion.
\end{proof}

\begin{corollary}\label{cor0} Let $N^{2n+1}(c)$ be a Sasakian space
form with constant $\varphi$-sectional curvature
$c\leq\frac{3-3n}{n+3}$. Then an integral submanifold $M^{n}$ with
constant mean curvature $|H|$ in $N^{2n+1}(c)$ is biharmonic if
and only if it is minimal.
\end{corollary}
\begin{proof} Assume that $M^{n}$ is a biharmonic integral submanifold with constant
mean curvature $|H|$ in $N^{2n+1}(c)$. It follows, from Theorem
\ref{t0}, that
$$
\begin{array}{ll}
g(\Delta^{\perp}H,H)&=-g(\trace B(\cdot,A_{H}\cdot),H)
+\frac{c(n+3)+3n-3}{4}|H|^{2}\\
\\&=\frac{c(n+3)+3n-3}{4}|H|^{2}-\sum_{i=1}^{n}g(B(X_{i},A_{H}X_{i}),H)\\
\\&=\frac{c(n+3)+3n-3}{4}|H|^{2}-\sum_{i=1}^{n}g(A_{H}X_{i},A_{H}X_{i})\\
\\&=\frac{c(n+3)+3n-3}{4}|H|^{2}-|A_{H}|^{2}.
\end{array}
$$
Thus, from the Weitzenb\"{o}ck formula
$$
\frac{1}{2}\Delta|H|^{2}=g(\Delta^{\perp}H,H)-|\nabla^{\perp}H|^{2},
$$
one obtains
\begin{equation}\label{eq:c}
\frac{c(n+3)+3n-3}{4}|H|^2-|A_H|^2-|\nabla^{\perp}H|^2=0.
\end{equation}
If $c<\frac{3-3n}{n+3}$, relation \eqref{eq:c} is equivalent to
$H=0$. Now, assume that $c=\frac{3-3n}{n+3}$. As for integral
submanifolds $\nabla^{\perp}H=0$ is equivalent to $H=0$, again
\eqref{eq:c} is equivalent to $H=0$.
\end{proof}

\begin{corollary}\label{cor0.1} Let $N^{2n+1}(c)$ be a Sasakian space
form with constant $\varphi$-sectional curvature
$c\leq\frac{3-3n}{n+3}$. Then a compact integral submanifold
$M^{n}$ is biharmonic if and only if it is minimal.
\end{corollary}
\begin{proof} Assume that $M^{n}$ is a biharmonic compact integral submanifold.
As in the proof of Corollary \ref{cor0} we have
$g(\Delta^{\perp}H,H)=\frac{c(n+3)+3n-3}{4}|H|^{2}-|A_{H}|^{2}$ and
so $\Delta|H|^{2}\leq 0$, which implies that $|H|^2=\cst$. Therefore
we obtain that $M$ is minimal in this case too.
\end{proof}

\begin{remark} From Corollary \ref{cor0} and Corollary
\ref{cor0.1} it is easy to see that in a Sasakian space form
$N^{2n+1}(c)$ with constant $\varphi$-sectional curvature $c\leq -3$
a biharmonic compact integral submanifold, or a biharmonic integral
submanifold $M^{n}$ with constant mean curvature, is minimal
whatever the dimension of $N$ is.
\end{remark}

\begin{proposition} Let $N^{2n+1}(c)$ be a Sasakian space
form and ${\bf i}:M^n\rightarrow N^{2n+1}$ be an integral
$\mathcal{C}$-parallel submanifold. Then $(\tau_2({\bf
i}))^{\top}=0$.
\end{proposition}
\begin{proof} Indeed, from Proposition \ref{prop0} we have
$|H|=\cst$ and $\nabla^{\perp}H\parallel\xi$, which implies that
$A_{\nabla^{\perp}_{X}H}=0$, for any vector field $X$ tangent to
$M$, since $A_{\xi}=0$, and so we conclude.
\end{proof}

\begin{proposition}\label{cor1} A non-minimal integral $\mathcal{C}$-parallel
submanifold $M^{n}$ of a Sasakian space form $N^{2n+1}(c)$ is
proper-biharmonic if and only if $c>\frac{7-3n}{n+3}$ and
$$
\trace B(\cdot,A_{H}\cdot)=\frac{c(n+3)+3n-7}{4}H.
$$
\end{proposition}
\begin{proof} We know, from Proposition \ref{prop0}, that $\Delta^{\perp}H=H$. Hence, from
Theorem \ref{t0} and the above Proposition, it follows that $M^n$
is biharmonic if and only if
$$
\trace B(\cdot,A_{H}\cdot)=\frac{c(n+3)+3n-7}{4}H.
$$

\noindent Next, if $M^n$ verifies the above condition, we contract
with $H$ and get
$$
|A_H|^2=\frac{c(n+3)+3n-7}{4}|H|^2.
$$
Since $A_H$ and $H$ do not vanish it follows that
$c>\frac{7-3n}{n+3}$.
\end{proof}

Now, let $\{X_i\}_{i=1}^n$ be an arbitrary orthonormal local frame
field on the integral $\mathcal{C}$-parallel submanifold $M^n$ of a
Sasakian space form $N^{2n+1}(c)$, and let $A_i=A_{\varphi X_i}$,
$i=\overline{1,n}$, be the corresponding shape operators. Then, from
Proposition \ref{cor1}, we obtain

\begin{proposition}\label{propmatrix} A non-minimal integral $\mathcal{C}$-parallel
submanifold $M^{n}$ of a Sasakian space form $N^{2n+1}(c)$,
$c>\frac{7-3n}{n+3}$, is proper-biharmonic if and only if
$$
\left(\begin{array}{cccc}g(A_1,A_1)&g(A_1,A_2)&\ldots&g(A_1,A_n)\\
\\g(A_2,A_1)&g(A_2,A_2)&\ldots&g(A_2,A_n)\\ \\\vdots&\vdots&\vdots&\vdots\\ \\
g(A_n,A_1)&g(A_n,A_2)&\ldots&g(A_n,A_n)\end{array}\right)\left(\begin{array}{c}
\trace A_1\\ \\\trace A_2\\ \\\vdots\\ \\\trace
A_n\end{array}\right)=k\left(\begin{array}{c} \trace A_1\\
\\\trace A_2\\ \\\vdots\\ \\\trace A_n\end{array}\right).
$$
where $k=\frac{c(n+3)+3n-7}{4}$.
\end{proposition}

\section{$3$-dimensional biharmonic integral $\mathcal{C}$-parallel
submanifolds of a Sasakian space form $N^{7}(c)$}
\setcounter{equation}{0}

In \cite{CBDBTK} Baikoussis, Blair and Koufogiorgios classified the
$3$-dimensional integral $\mathcal{C}$-parallel submanifolds in a
Sasakian space form $(N^{7}(c),\varphi,\xi,\eta,g)$. In order to
obtain the classification, they worked with a special local
orthonormal basis (see also \cite{FDLV}). Here we shall briefly
recall how this basis is constructed.

Let ${\bf i}:M^{3}\rightarrow N^{7}(c)$ be an integral submanifold
of non-zero constant mean curvature. Let $p$ be an arbitrary point
of $M$, and consider the function $f_p:U_{p}M\rightarrow
\mathbb{R}$ given by
$$
f_p(u)=g(B(u,u),\varphi u),
$$
where $U_{p}M=\{u\in T_{p}M:g(u,u)=1\}$ is the unit sphere in the
tangent space $T_{p}M$. If $f_p(u)=0$, for all $u\in U_pM$, then,
for any $v_1, v_2\in U_pM$ such that $g(v_1,v_2)=0$ we have that
$$
g(B(v_1,v_1),\varphi v_1)=0,\quad g(B(v_1,v_1),\varphi
v_2)=0,\quad g(B(v_1,v_2),\varphi v_1)=0.
$$
Now, if $\{X_1,X_2,X_3\}$ is an arbitrary orthonormal basis at
$p$, it follows that $\trace A_{\varphi X_i}$ $=0$, for any
$i=\overline{1,3}$, and therefore $H(p)=0$. Consequently, the
function $f_p$ does not vanish identically.

Since $U_{p}M$ is compact, $f_p$ attains an absolute maximum at a
unit vector $X_{1}$. It follows that
$$
\begin{cases}
g(B(X_1,X_1),\varphi X_1)>0,\quad g(B(X_1,X_1),\varphi X_1)\geq
|g(B(w,w),\varphi w)|\\
g(B(X_1,X_1),\varphi w)=0, \quad g(B(X_1,X_1),\varphi X_1)\geq
2g(B(w,w),\varphi X_1),
\end{cases}
$$
where $w$ is a unit vector tangent to $M$ at $p$ and orthogonal to
$X_1$. It is easy to see that $X_1$ is an eigenvector of
$A_{1}=A_{\varphi X_{1}}$ with corresponding eigenvalue
$\lambda_1$. Then, since $A_{1}$ is symmetric, we consider $X_{2}$
and $X_{3}$ to be unit eigenvectors of $A_{1}$ orthogonal to each
other and to $X_1$. Further, we distinguish two cases.

If $\lambda_2\neq\lambda_3$, we can choose $X_2$ and $X_3$ such
that
$$
\begin{cases}
g(B(X_2,X_2),\varphi X_2)\geq 0,\quad g(B(X_3,X_3),\varphi
X_3)\geq 0\\ g(B(X_2,X_2),\varphi X_2)\geq g(B(X_3,X_3),\varphi
X_3).
\end{cases}
$$

If $\lambda_2=\lambda_3$, we consider $f_{1,p}$ the restriction of
$f_p$ to $\{w\in U_{p}M: g(w,X_{1})=0\}$, and we have two
subcases:
\begin{enumerate}

\item the function $f_{1,p}$ is identically zero. In this case, we
have
$$
\begin{cases}
g(B(X_2,X_2),\varphi X_2)=0,\quad g(B(X_2,X_2),\varphi X_3)=0\\
g(B(X_2,X_3),\varphi X_3)=0,\quad g(B(X_3,X_3),\varphi X_3)=0.
\end{cases}
$$

\item the function $f_{1,p}$ does not vanish identically. Then we
choose $X_2$ such that $f_{1,p}(X_2)$ is an absolute maximum. We
have that
$$
\begin{cases}
g(B(X_2,X_2),\varphi X_2)>0, \quad g(B(X_2,X_2),\varphi X_2)\geq
g(B(X_3,X_3),\varphi X_3)\geq 0
\\
g(B(X_2,X_2),\varphi X_3)=0, \quad g(B(X_2,X_2),\varphi X_2)\geq
2g(B(X_3,X_3),\varphi X_2).
\end{cases}
$$
\end{enumerate}
Now, with respect to the orthonormal basis
$\{X_{1},X_{2},X_{3}\}$, the shape operators $A_{1}$,
$A_{2}=A_{\varphi X_{2}}$ and $A_{3}=A_{\varphi X_{3}}$, at $p$,
can be written as follows
\begin{equation}\label{eq:4.1}
A_{1}=\left(\begin{array}{ccc}\lambda_{1}&0&0\\
0&\lambda_{2}&0\\
0&0&\lambda_{3}\end{array}\right),\quad A_{2}=\left(\begin{array}{ccc}0&\lambda_{2}&0\\
\lambda_{2}&\alpha&\beta\\
0&\beta&\gamma\end{array}\right),\quad A_{3}=\left(\begin{array}{ccc}0&0&\lambda_{3}\\
0&\beta&\gamma\\
\lambda_{3}&\gamma&\delta\end{array}\right).
\end{equation}
We also have $A_0=A_{\xi}=0$. With these notations we have
\begin{equation}\label{eq:basis1}
\lambda_1>0,\quad \lambda_1\geq |\alpha|,\quad \lambda_1\geq
|\delta|,\quad \lambda_1\geq 2\lambda_2,\quad \lambda_1\geq
2\lambda_3.
\end{equation}

For $\lambda_2\neq\lambda_3$ we get
\begin{equation}\label{eq:basis2}
\alpha\geq 0,\quad\delta\geq 0\quad\hbox{and}\quad
\alpha\geq\delta
\end{equation}
and for $\lambda_2=\lambda_3$ we obtain that
\begin{equation}
\alpha=\beta=\gamma=\delta=0
\end{equation}
or
\begin{equation}\label{eq:basis3}
\alpha>0,\quad\delta\geq 0,\quad\alpha\geq\delta,\quad
\beta=0\quad\hbox{and}\quad\alpha\geq 2\gamma.
\end{equation}

We can extend $X_1$ on a neighbourhood $V_p$ of $p$ such that
$X_1(q)$ is a maximal point of $f_q:U_qM\rightarrow \mathbb{R}$,
for any point $q$ of $V_p$.

If the eigenvalues of $A_{1}$ have constant multiplicities, then the
above basis $\{X_{1},X_{2},$ $X_{3}\}$, defined at $p$, can be
smoothly extended and we can work on the open dense subset of $M$
defined by this property.

Using this basis, in \cite{CBDBTK}, the authors proved that, when
$M$ is an integral $\mathcal{C}$-parallel submanifold, the functions
$\lambda_i$, $i=\overline{1,3}$, and $\alpha$, $\beta$, $\gamma$,
$\delta$ are constant on $V_p$, and then classified all
$3$-dimensional integral $\mathcal{C}$-parallel submanifolds in a
$7$-dimensional Sasakian space form.

According to that classification, if $c>-3$ then $M$ is a
non-minimal integral $\mathcal{C}$-parallel submanifold if and only
if either:

\noindent {\bf Case I.} $M$ is flat, locally it is a product of
three curves, which are helices of osculating orders $r\leq 4$,
and $\lambda_{1}=\frac{\lambda^{2}-\frac{c+3}{4}}{\lambda}$,
$\lambda_{2}=\lambda_{3}=\lambda=\cst\neq 0$, $\alpha=\cst$,
$\beta=0$, $\gamma=\cst$, $\delta=\cst$, such that
$-\frac{\sqrt{c+3}}{2}<\lambda< 0$, $0<\alpha\leq\lambda_1$,
$\alpha>2\gamma$, $\alpha\geq\delta\geq 0$,
$\frac{c+3}{4}+\lambda^{2}+\alpha\gamma-\gamma^{2}=0$ and
$\Big(\frac{3\lambda^2-\frac{c+3}{4}}{\lambda}\Big)^2+(\alpha+\gamma)^2+\delta^2>0$.

\noindent {\bf Case II.} $M$ is locally isometric to a product
$\Gamma\times\bar{M}^{2}$, where $\Gamma$ is a curve and
$\bar{M}^{2}$ is a $\mathcal{C}$-parallel surface, and either
\begin{itemize}
\item
[(1)]$\lambda_{1}=2\lambda_{2}=\frac{\sqrt{c+3}}{2\sqrt{2}}$,
$\lambda_{3}=-\frac{\sqrt{c+3}}{2\sqrt{2}}$,
$\alpha=\gamma=\delta=0$,
$\beta=\pm\frac{\sqrt{3(c+3)}}{4\sqrt{2}}$. In this case $\Gamma$ is
a helix in $N$ with curvatures $\kappa_{1}=\frac{1}{\sqrt{2}}$ and
$\kappa_{2}=1$, and $\bar{M}^{2}$ is locally isometric to the
$2$-dimensional Euclidean sphere of radius
$\rho=\sqrt{\frac{8}{3(c+3)}}$.

\noindent or \item
[(2)]$\lambda_{1}=\frac{\lambda^{2}-\frac{c+3}{4}}{\lambda}$,
$\lambda_{2}=\lambda_{3}=\lambda=\cst$,
$\alpha=\beta=\gamma=\delta=0$, such that
$-\frac{\sqrt{c+3}}{2}<\lambda< 0$ and
$\lambda^2\neq\frac{c+3}{12}$. In this case $\Gamma$ is a helix in
$N$ with curvatures $\kappa_{1}=\lambda_{1}$ and $\kappa_{2}=1$, and
$\bar{M}^{2}$ is the $2$-dimensional Euclidean sphere of radius
$\rho=\frac{1}{\sqrt{\frac{c+3}{4}+\lambda^{2}}}$.
\end{itemize}

Now, identifying the shape operators $A_i$ with the corresponding
matrices, from Proposition \ref{propmatrix}, we get

\begin{proposition}\label{propmatrix37} A non-minimal integral $\mathcal{C}$-parallel
submanifold $M^3$ of a Sasakian space form $N^7(c)$,
$c>-\frac{1}{3}$, is proper-biharmonic if and only if
\begin{equation}\label{eq:matrix}
\Bigg(\sum_{i=1}^3A_i^2\Bigg)\left(\begin{array}{c}\trace
A_1\\\trace A_2\\\trace
A_3\end{array}\right)=\frac{3c+1}{2}\left(\begin{array}{c}\trace
A_1\\\trace A_2\\\trace A_3\end{array}\right),
\end{equation}
where matrices $A_i$ are given by \eqref{eq:4.1}.
\end{proposition}

Now, we can state

\begin{theorem}\label{t1} A $3$-dimensional integral $\mathcal{C}$-parallel submanifold
$M^{3}$ of a Sasakian space form $N^{7}(c)$ is proper-biharmonic
if and only if either
\begin{enumerate}
\item $c>-\frac{1}{3}$ and $M^{3}$ is flat and locally is a
product of three curves:
\begin{itemize}
\item The $X_{1}$-curve is a helix with curvatures
$\kappa_{1}=\frac{\lambda^{2}-\frac{c+3}{4}}{\lambda}$ and
$\kappa_{2}=1$,

\item The $X_{2}$-curve is a helix of order $4$ with curvatures
$\kappa_{1}=\sqrt{\lambda^{2}+\alpha^{2}}$,
$\kappa_{2}=\frac{\alpha}{\kappa_{1}}\sqrt{\lambda^{2}+1}$ and
$\kappa_{3}=-\frac{\lambda\sqrt{\lambda^{2}+1}}{\kappa_{1}}$,

\item The $X_{3}$-curve is a helix of order $4$ with curvatures
$\kappa_{1}=\sqrt{\lambda^{2}+\gamma^{2}+\delta^{2}}$,
$\kappa_{2}=\frac{\delta}{\kappa_{1}}\sqrt{\lambda^{2}+\gamma^{2}+1}$
and
$\kappa_{3}=\frac{\kappa_{2}\sqrt{\lambda^{2}+\gamma^{2}}}{\delta}$,
if $\delta\neq 0$, or a circle with curvature
$\kappa_{1}=\sqrt{\lambda^{2}+\gamma^{2}}$, if $\delta=0$,
\end{itemize}
\noindent where $\lambda,\alpha,\gamma,\delta$ are constants given
by
\begin{equation}\label{system0}
\begin{cases}
(3\lambda^{2}-\frac{c+3}{4})\Big(3\lambda^{4}-2(c+1)\lambda^{2}+
\frac{(c+3)^{2}}{16}\Big)+\lambda^{4}((\alpha+\gamma)^{2}+\delta^{2})=0\\
(\alpha+\gamma)(5\lambda^{2}+\alpha^{2}+\gamma^{2}-\frac{7c+5}{4})+\gamma\delta^{2}=0\\
\delta(5\lambda^{2}+\delta^{2}+3\gamma^{2}+\alpha\gamma-\frac{7c+5}{4})=0\\
\frac{c+3}{4}+\lambda^{2}+\alpha\gamma-\gamma^{2}=0
\end{cases}
\end{equation}
such that $-\frac{\sqrt{c+3}}{2}<\lambda< 0$,
$0<\alpha\leq\frac{\lambda^2-\frac{c+3}{4}}{\lambda}$,
$\alpha\geq\delta\geq 0$, $\alpha>2\gamma$ and
$\lambda^{2}\neq\frac{c+3}{12}$;

\noindent or

\item $M^{3}$ is locally isometric to a product
$\Gamma\times\bar{M}^{2}$, between a curve and a
$\mathcal{C}$-parallel surface of $N$, and either:
\begin{enumerate}
\item $c=\frac{5}{9}$, $\Gamma$ is a helix in
$N^{7}(\frac{5}{9})$ with curvatures
$\kappa_{1}=\frac{1}{\sqrt{2}}$ and $\kappa_{2}=1$, and
$\bar{M}^{2}$ is locally isometric to the $2$-dimensional
Euclidean sphere with radius $\frac{\sqrt{3}}{2}$;

\noindent or

\item
$c\in\Big[\frac{-7+8\sqrt{3}}{13},+\infty\Big)\setminus\{1\}$,
$\Gamma$ is a helix in $N^{7}(c)$ with curvatures
$\kappa_{1}=\frac{\lambda^{2}-\frac{c+3}{4}}{\lambda}$ and
$\kappa_{2}=1$, and $\bar{M}^{2}$ is locally isometric to the
$2$-dimensional Euclidean sphere with radius
$\frac{2}{\sqrt{4\lambda^{2}+c+3}}$, where
\begin{equation}\label{eq:4.4}
\lambda^{2}=\begin{cases}\frac{4c+4\pm\sqrt{13c^{2}+14c-11}}{12}
&\text{if}\quad c<1\\ \frac{4c+4-\sqrt{13c^{2}+14c-11}}{12}
&\text{if}\quad c>1\end{cases}\quad\text{and}\quad\lambda<0.
\end{equation}
\end{enumerate}
\end{enumerate}
\end{theorem}
\begin{proof} Let $M^{3}$ be a proper-biharmonic integral $\mathcal{C}$-parallel
submanifold of a Sasakian space form $N^{7}(c)$. From Proposition
\ref{propmatrix37} we see that $c>-\frac{1}{3}$.

Next, we easily get that the equation \eqref{eq:matrix} is
equivalent to the system
\begin{equation}\label{system}
\begin{cases}
(\sum_{i=1}^{3}\lambda_{i})(\sum_{i=1}^{3}\lambda_{i}^{2}-\frac{3c+1}{2})+(\alpha+\gamma)(\alpha\lambda_{2}+
\gamma\lambda_{3})\\+(\beta+\delta)(\beta\lambda_{2}+\delta\lambda_{3})=0\\
(\sum_{i=1}^{3}\lambda_{i})(\alpha\lambda_{2}+
\gamma\lambda_{3})+(\alpha+\gamma)(2\lambda_{2}^{2}+\alpha^{2}+3\beta^{2}+\gamma^{2}+
\beta\delta-\frac{3c+1}{2})\\+\gamma(\beta+\delta)^{2}=0\\
(\sum_{i=1}^{3}\lambda_{i})(\beta\lambda_{2}+
\delta\lambda_{3})+\beta(\alpha+\gamma)^{2}\\
+(\beta+\delta)(2\lambda_{3}^{2}+\delta^{2}+3\gamma^{2}+\beta^{2}+\alpha\gamma-\frac{3c+1}{2})=0.
\end{cases}
\end{equation}

In the following, we shall split the study of this system, as $M^3$
is given by \textbf{Case I} or \textbf{Case II} of the
classification.

\noindent \textbf{Case I.} The system \eqref{system} is equivalent
to the system given by the first three equations of \eqref{system0}.
Now, $M$ is not minimal if and only if at least one of the
components of the mean curvature vector field $H$ does not vanish
and, from the first equation of \eqref{system0}, it follows that
$\lambda^2$ must be different from $\frac{c+3}{12}$. Thus, again
using \cite{CBDBTK} for the expressions of the curvatures of the
three curves, we obtain the first case of the Theorem.

\noindent \textbf{Case II.}
\begin{enumerate}
\item[(1)] The first and the third equation of \eqref{system} are
equivalent, in this case, to $c=\frac{5}{9}$ and the second equation
is identically satisfied. Then, from the classification of the
integral $\mathcal{C}$-parallel submanifolds, we get the first part
of the second case of the Theorem.

\item[(2)] The second and the third equation of system
\eqref{system} are satisfied, in this case, and the first equation
is equivalent to
$$
3\lambda^{4}-2(c+1)\lambda^{2}+\frac{(c+3)^{2}}{16}=0.
$$
This equation has solutions if and only if
$$
c\in\Big(-\infty,\frac{-7-8\sqrt{3}}{13}\Big]\cup\Big[\frac{-7+8\sqrt{3}}{13},+\infty\Big),
$$
and these solutions are given by
$$
\lambda^{2}=\frac{4c+4\pm\sqrt{13c^{2}+14c-11}}{12}.
$$
Since $c>-\frac{1}{3}$ it follows that
$c\in\Big[\frac{-7+8\sqrt{3}}{13},+\infty\Big)$. Moreover, if $c=1$,
from the above relation, it follows that $\lambda^2$ must be equal
to $1$ or $\frac{1}{3}$, which is a contradiction, and therefore
$c\in\Big[\frac{-7+8\sqrt{3}}{13},+\infty\Big)\setminus\{1\}$.
Further, it is easy to check that
$\lambda^{2}=\frac{4c+4+\sqrt{13c^{2}+14c-11}}{12}<\frac{c+3}{4}$ if
and only if $c\in\Big[\frac{-7+8\sqrt{3}}{13},1)$ and
$\lambda^{2}=\frac{4c+4-\sqrt{13c^{2}+14c-11}}{12}<\frac{c+3}{4}$ if
and only if
$c\in\Big[\frac{-7+8\sqrt{3}}{13},+\infty\Big)\setminus\{1\}$.
\end{enumerate}
\end{proof}

\section{Proper-biharmonic submanifolds in the 7-sphere}
\setcounter{equation}{0}

In this section we shall work with the standard model for simply
connected Sasakian space forms $N^{7}(c)$ with $c>-3$, which is the
sphere $\mathbb{S}^{7}$ endowed with its canonical Sasakian
structure or with the deformed Sasakian structure introduced by
Tanno.

In \cite{CBDBTK} the authors obtained the explicit equation of the
$3$-dimensional integral $\mathcal{C}$-parallel flat submanifolds in
$\mathbb{S}^{7}(1)$, whilst in \cite{DFCO2} we gave the explicit
equation of such submanifolds in $\mathbb{S}^{7}(c)$, $c>-3$.

Using these results and Theorem \ref{t1} we easily get

\begin{theorem}\label{t2} A $3$-dimensional integral $\mathcal{C}$-parallel submanifold
$M^{3}$ of $\mathbb{S}^{7}(c)$, $c=\frac{4}{a}-3>-3$, is
proper-biharmonic if and only if either
\begin{enumerate}
\item $c>-\frac{1}{3}$ and $M^{3}$ is flat, locally is a product
of three curves and its position vector in $\mathbb{C}^{4}$ is
$$
\begin{array}{ll}
x(u,v,w)=&\frac{\lambda}{\sqrt{\lambda^{2}+\frac{1}{a}}}\exp(\mathrm{i}(\frac{1}{a\lambda}u))\mathcal{E}_{1}
+\frac{1}{\sqrt{a(\gamma-\alpha)(2\gamma-\alpha)}}\exp(-\mathrm{i}(\lambda
u-(\gamma-\alpha)v))\mathcal{E}_{2}\\
\\&+\frac{1}{\sqrt{a\rho_{1}(\rho_{1}+\rho_{2})}}\exp(-\mathrm{i}(\lambda u+\gamma
v+\rho_{1}w))\mathcal{E}_{3}\\
\\&+\frac{1}{\sqrt{a\rho_{2}(\rho_{1}+\rho_{2})}}\exp(-\mathrm{i}
(\lambda u+\gamma v-\rho_{2}w))\mathcal{E}_{4},
\end{array}
$$
where
$\rho_{1,2}=\frac{1}{2}(\sqrt{4\gamma(2\gamma-\alpha)+\delta^{2}}\pm\delta)$
and $\lambda,\alpha,\gamma,\delta$ are real constants given by
\eqref{system0} such that $-\frac{1}{\sqrt{a}}<\lambda<0$,
$0<\alpha\leq\frac{\lambda^2-\frac{1}{a}}{\lambda}$,
$\alpha\geq\delta\geq 0$, $\alpha>2\gamma$,
$\lambda^{2}\neq\frac{1}{3a}$ and $\{\mathcal{E}_i\}_{i=1}^4$ is an
orthonormal basis of $\mathbb{C}^4$ with respect to the usual
Hermitian inner product;

\noindent or

\item $M^{3}$ is locally isometric to a product
$\Gamma\times\bar{M}^{2}$, between a curve and a
$\mathcal{C}$-parallel surface of $N$, and either:
\begin{enumerate}
\item $c=\frac{5}{9}$, $\Gamma$ is a helix in
$\mathbb{S}^{7}\Big(\frac{5}{9}\Big)$ with curvatures
$\kappa_{1}=\frac{1}{\sqrt{2}}$ and $\kappa_{2}=1$, and
$\bar{M}^{2}$ is locally isometric to the $2$-dimensional
Euclidean sphere with radius $\frac{\sqrt{3}}{2}$;

\noindent or

\item
$c\in\Big[\frac{-7+8\sqrt{3}}{13},+\infty\Big)\setminus\{1\}$,
$\Gamma$ is a helix in $\mathbb{S}^{7}(c)$ with curvatures
$\kappa_{1}=\frac{\lambda^{2}-\frac{c+3}{4}}{\lambda}$ and
$\kappa_{2}=1$, and $\bar{M}^{2}$ is locally isometric to the
$2$-dimensional Euclidean sphere with radius
$\frac{2}{\sqrt{4\lambda^{2}+c+3}}$, where
$$
\lambda^{2}=\begin{cases}\frac{4c+4\pm\sqrt{13c^{2}+14c-11}}{12}
&\text{if}\quad  c<1\\
\frac{4c+4-\sqrt{13c^{2}+14c-11}}{12}&\text{if}\quad
c>1\end{cases}\quad\text{and}\quad\lambda<0.
$$
\end{enumerate}
\end{enumerate}
\end{theorem}

Now, applying this Theorem in the case of the $7$-sphere endowed
with its canonical Sasakian structure we get the following
Corollary, which also shows that, for $c=1$, the system
\eqref{system0} can be completely solved.

\begin{corollary}
\label{corolarfinal} A $3$-dimensional integral
$\mathcal{C}$-parallel submanifold $M^{3}$ of $\mathbb{S}^{7}(1)$ is
proper-biharmonic if and only if it is flat, locally it is a product
of three curves and its position vector in $\mathbb{C}^{4}$ is
$$
\begin{array}{ll}
x(u,v,w)=&-\frac{1}{\sqrt{6}}\exp(-\mathrm{i}\sqrt{5}u)\mathcal{E}_{1}
+\frac{1}{\sqrt{6}}\exp(\mathrm{i}(\frac{1}{\sqrt{5}}
u-\frac{4\sqrt{3}}{\sqrt{10}}v))\mathcal{E}_{2}\\
\\&+\frac{1}{\sqrt{6}}\exp(\mathrm{i}(\frac{1}{\sqrt{5}}u+\frac{\sqrt{3}}{\sqrt{10}}
v-\frac{3\sqrt{2}}{2}w))\mathcal{E}_{3}\\
\\&+\frac{1}{\sqrt{2}}\exp(\mathrm{i}(\frac{1}{\sqrt{5}}u+\frac{\sqrt{3}}{\sqrt{10}}
v+\frac{\sqrt{2}}{2}w))\mathcal{E}_{4},
\end{array}
$$
where $\{\mathcal{E}_i\}_{i=1}^4$ is an orthonormal basis of
$\mathbb{C}^4$ with respect to the usual Hermitian inner product.
Moreover, the $X_1(=x_u)$-curve is a helix with curvatures
$\kappa_{1}=\frac{4\sqrt{5}}{5}$ and $\kappa_{2}=1$, the
$X_2(=x_v)$-curve is a helix of order $4$ with curvatures
$\kappa_{1}=\frac{\sqrt{29}}{\sqrt{10}}$,
$\kappa_{2}=\frac{9\sqrt{2}}{\sqrt{145}}$ and
$\kappa_{3}=\frac{2\sqrt{3}}{\sqrt{145}}$ and the
$X_3(=x_w)$-curve is a helix of order $4$ with curvatures
$\kappa_{1}=\frac{\sqrt{5}}{\sqrt{2}}$,
$\kappa_{2}=\frac{2\sqrt{3}}{\sqrt{10}}$ and
$\kappa_{3}=\frac{\sqrt{3}}{\sqrt{10}}$.
\end{corollary}

\begin{proof} Since $c=1$ the system \eqref{system0} becomes
\begin{equation}\label{system1}
\begin{cases}
(3\lambda^{2}-1)^2(\lambda^2-1)+\lambda^{4}((\alpha+\gamma)^{2}+\delta^{2})=0\\
(\alpha+\gamma)(5\lambda^{2}+\alpha^{2}+\gamma^{2}-3)+\gamma\delta^{2}=0\\
\delta(5\lambda^{2}+\delta^{2}+3\gamma^{2}+\alpha\gamma-3)=0\\
\lambda^{2}+\alpha\gamma-\gamma^{2}+1=0
\end{cases}
\end{equation}
with the supplementary conditions
\begin{equation}\label{conditions}
-1<\lambda<0,\quad 0<\alpha\leq\frac{\lambda^2-1}{\lambda},\quad
\alpha\geq\delta\geq 0,\quad\alpha>2\gamma\quad\text{and}\quad
\lambda^{2}\neq\frac{1}{3}.
\end{equation}
We note that, since $\alpha>2\gamma$, from the fourth equation of
\eqref{system1} it results that $\gamma<0$.

The third equation of system \eqref{system1} suggests that, in
order to solve this system, we need to split our study in two
cases as $\delta$ is equal to $0$ or not.

\noindent\textbf{Case 1: $\delta=0$.} In this case the third
equation holds whatever the values of $\lambda$, $\alpha$ and
$\gamma$ are, and so does the condition $\alpha\geq\delta$. We
also note that $\alpha\neq-\gamma$, since otherwise, from the
first equation, it results $\lambda^2=1$ or
$\lambda^2=\frac{1}{3}$, which are both contradictions.

In the following, we shall look for $\alpha$ of the form
$\alpha=\omega\gamma$, where
$\omega\in(-\infty,0)\setminus\{-1\}$. From the second and the
fourth equations of the system we have
$\lambda^2=-\frac{\omega^2+3\omega-2}{(\omega-2)(\omega-3)}$
$\gamma^2=\frac{8}{(\omega-2)(\omega-3)}$ and then
$\alpha^2=\frac{8\omega^2}{(\omega-2)(\omega-3)}$. Replacing in
the first equation, after a straightforward computation, it can be
written as
$$
\frac{8(\omega+1)^3(1-3\omega)}{(\omega-3)^3(\omega-2)}=0
$$
and its solutions are $-1$ and $\frac{1}{3}$. But
$\omega\in(-\infty,0)\setminus\{-1\}$ and therefore we conclude
that there are no solutions of the system that verify all
conditions \eqref{conditions} when $\delta=0$.

\noindent\textbf{Case 2: $\delta>0$.} In this case the third
equation of \eqref{system1} becomes
$$
5\lambda^2+\delta^2+3\gamma^2+\alpha\gamma-3=0.
$$
Now, again taking $\alpha=\omega\gamma$, this time with
$\omega\in(-\infty,0)$, from the last three equations of the
system, we easily get
$\lambda^2=-\frac{\omega^2+5\omega+2}{(\omega-1)(\omega-2)}$,
$\alpha^2=\frac{8\omega^3}{(\omega-1)^2(\omega-2)}$,
$\gamma^2=\frac{8\omega}{(\omega-1)^2(\omega-2)}$ and
$\delta^2=\frac{8(\omega+1)^2}{(\omega-1)^2}$.

\noindent Next, from the first equation of \eqref{system1}, after
a straightforward computation, one obtains
$$
\frac{16(\omega+1)^3(\omega+3)}{(\omega-2)(\omega-1)^3}=0,
$$
which solutions are $-3$ and $-1$. If $\omega=-1$ it follows that
$\lambda^2=\frac{1}{3}$, which is a contradiction, and therefore
we obtain that $\omega=-3$. Hence
$$
\lambda^2=\frac{1}{5},\quad\alpha^2=\frac{27}{10},\quad\gamma^2=\frac{3}{10}\quad\text{and}\quad\delta^2=2.
$$
As $\lambda<0$, $\alpha>0$, $\gamma<0$ and $\delta>0$ it results
that $\lambda=-\frac{1}{\sqrt{5}}$,
$\alpha=\frac{3\sqrt{3}}{\sqrt{10}}$,
$\gamma=-\frac{\sqrt{3}}{\sqrt{10}}$ and $\delta=\sqrt{2}$. It can
be easily seen that also the conditions \eqref{conditions} are
verified by these values, and then, by the meaning of the first
statement of Theorem \ref{t2}, we come to the conclusion.
\end{proof}

\begin{remark} A proper-biharmonic compact submanifold $M$ of
$\mathbb{S}^n$ of constant mean curvature $|H|\in(0,1)$ is of
$2$-type and mass-symmetric (see~\cite{ABSMCO,ABSMCO2}). In our
case, the Riemannian immersion $x$ can be written as $x=x_1+x_2$,
where
$$
x_1(u,v,w)=\frac{1}{\sqrt{2}}\exp(\mathrm{i}(\frac{1}{\sqrt{5}}u+\frac{\sqrt{3}}{\sqrt{10}}
v+\frac{\sqrt{2}}{2}w))\mathcal{E}_{4},
$$
$$
\begin{array}{ll}
x_2(u,v,w)=&-\frac{1}{\sqrt{6}}\exp(-\mathrm{i}\sqrt{5}u)\mathcal{E}_{1}
+\frac{1}{\sqrt{6}}\exp(\mathrm{i}(\frac{1}{\sqrt{5}}
u-\frac{4\sqrt{3}}{\sqrt{10}}v))\mathcal{E}_{2}\\
\\&+\frac{1}{\sqrt{6}}\exp(\mathrm{i}(\frac{1}{\sqrt{5}}u+\frac{\sqrt{3}}{\sqrt{10}}
v-\frac{3\sqrt{2}}{2}w))\mathcal{E}_{3},\\
\end{array}
$$
and $\Delta x_1=3(1-|H|)x_1=x_1$, $\Delta x_2=3(1+|H|)x_2=5x_2$,
$|H|=\frac{2}{3}$. Now, Corollary~\ref{corolarfinal} could also be
proved by using the main result in ~\cite{CBDB2} and Proposition 4.1
in ~\cite{RCSMCO}.
\end{remark}

\begin{remark} By a straightforward computation we can deduce that
the map $x$ factorizes to a map from the torus
$\mathcal{T}^3=\mathbb{R}^3/\Lambda$ into $\mathbb{R}^8$, where
$\Lambda$ is the lattice generated by the vectors
$a_1=(\frac{6\pi}{\sqrt{5}},\frac{\sqrt{3}\pi}{\sqrt{10}},\frac{\pi}{\sqrt{2}})$,
$a_2=(0,-\frac{3\sqrt{5}\pi}{\sqrt{6}},-\frac{\pi}{\sqrt{2}})$ and
$a_3=(0,0,-\frac{4\pi}{\sqrt{2}})$, and the quotient map is a
Riemannian immersion.
\end{remark}

By the meaning of Theorem \ref{teorema1} we know that the cylinder
over $x$, given by
$$
y(t,u,v,w)=\phi_t(x(u,v,w)),
$$
is a proper-biharmonic map into $\mathbb{S}^7(1)$. Moreover, we
have

\begin{proposition} The cylinder over $x$ determines a
proper-biharmonic Riemannian embedding from the torus
$\mathcal{T}^4=\mathbb{R}^4/\Lambda$ into $\mathbb{S}^7$, where the
lattice $\Lambda$ is generated by
$a_1=(\frac{2\pi}{\sqrt{6}},0,0,0)$
$a_2=(0,\frac{2\pi}{\sqrt{6}},0,0)$,
$a_3=(0,0,\frac{2\pi}{\sqrt{6}},0)$ and
$a_4=(0,0,0,\frac{2\pi}{\sqrt{2}})$. The image of this embedding is
the Riemannian product between a Euclidean circle of radius
$\frac{1}{\sqrt{2}}$ and three other Euclidean circles, each of
radius $\frac{1}{\sqrt{6}}$.
\end{proposition}

\begin{proof} As the flow of the characteristic vector field $\xi$ is given
by $\phi_t(z)=\exp(-\mathrm{i}t)z$ we get
$$
\begin{array}{ll}
y(t,u,v,w)=&-\frac{1}{\sqrt{6}}\exp(-\mathrm{i}(t+\sqrt{5}u))\mathcal{E}_{1}
+\frac{1}{\sqrt{6}}\exp(\mathrm{i}(-t+\frac{1}{\sqrt{5}}
u-\frac{4\sqrt{3}}{\sqrt{10}}v))\mathcal{E}_{2}\\
\\&+\frac{1}{\sqrt{6}}\exp(\mathrm{i}(-t+\frac{1}{\sqrt{5}}u+\frac{\sqrt{3}}{\sqrt{10}}
v-\frac{3\sqrt{2}}{2}w))\mathcal{E}_{3}\\
\\&+\frac{1}{\sqrt{2}}\exp(\mathrm{i}(-t+\frac{1}{\sqrt{5}}u+\frac{\sqrt{3}}{\sqrt{10}}
v+\frac{\sqrt{2}}{2}w))\mathcal{E}_{4},
\end{array}
$$
where $\{\mathcal{E}_i\}_{i=1}^4$ is an orthonormal basis of
$\mathbb{C}^4$ with respect to the usual Hermitian inner product.

Now, we consider the following two orthogonal transformations of
$\mathbb{R}^4$:
$$
\begin{cases}
\frac{1}{\sqrt{2}}t+\frac{1}{\sqrt{10}}u+\frac{\sqrt{3}}{2\sqrt{5}}
v+\frac{1}{2}w=t'\\\frac{2}{\sqrt{5}}u-\frac{\sqrt{6}}{4\sqrt{5}}v-\frac{\sqrt{2}}{4}w=u'\\
\frac{\sqrt{5}}{2\sqrt{2}}v-\frac{\sqrt{3}}{2\sqrt{2}}w=v'\\
\frac{1}{\sqrt{2}}t-\frac{1}{\sqrt{10}}u-\frac{\sqrt{3}}{2\sqrt{5}}
v-\frac{1}{2}w=w'
\end{cases}
\quad\text{and}\quad
\begin{cases}
\frac{\sqrt{2}}{\sqrt{6}}t'+\frac{2}{\sqrt{6}}u'=\widetilde{t}\\
-\frac{\sqrt{2}}{\sqrt{6}}t'+\frac{1}{\sqrt{6}}u'-\frac{\sqrt{3}}{\sqrt{6}}v'=\widetilde{u}\\
-\frac{\sqrt{2}}{\sqrt{6}}t'+\frac{1}{\sqrt{6}}u'+\frac{\sqrt{3}}{\sqrt{6}}v'=\widetilde{v}\\
w'=\widetilde{w}
\end{cases}
$$
and obtain
$$
\begin{array}{ll}
\widetilde{y}(\widetilde{t},\widetilde{u},\widetilde{v},\widetilde{w})
=&-\frac{1}{\sqrt{6}}\exp(-\mathrm{i}(\sqrt{6}\widetilde{t}))\mathcal{E}_{1}
+\frac{1}{\sqrt{6}}\exp(\mathrm{i}(\sqrt{6}\widetilde{u}))\mathcal{E}_{2}
+\frac{1}{\sqrt{6}}\exp(\mathrm{i}(\sqrt{6}\widetilde{v}))\mathcal{E}_{3}\\ \\
&+\frac{1}{\sqrt{2}}\exp(\mathrm{i}(\sqrt{2}\widetilde{w}))\mathcal{E}_{4},
\end{array}
$$
which ends the proof.
\end{proof}

\begin{remark}
We see that $y$ can be written as $y=y_1+y_2$, where
$y_1(t,u,v,w)=\exp(-\mathrm{i}t)x_1$,
$y_2(t,u,v,w)=\exp(-\mathrm{i}t)x_2$, and $\Delta y_1=2y_1$, $\Delta
y_2=6y_2$, the mean curvature of $y$ being equal to $\frac{1}{2}$.
\end{remark}

\begin{remark} It is known that the parallel flat
$(n+1)$-dimensional compact anti-invariant submanifolds in
$\mathbb{S}^{2n+1}(1)$ are Riemannian products of circles of radii
$r_i$, $i=\overline{1,n+1}$, where $\sum_{i=1}^{n+1}r_i^2=1$ (see
\cite{KYMK}). The biharmonicity of such submanifolds was solved in
\cite{WZ}.
\end{remark}

\section{Proper-biharmonic parallel Lagrangian submanifolds of
$\mathbb{C}P^3$}

We consider the Hopf fibration
$\pi:\mathbb{S}^{2n+1}(1)\to\mathbb{C}P^n(4)$, and $\overline{M}$ a
Lagrangian submanifold of $\mathbb{C}P^n$. Then
$\widetilde{M}=\pi^{-1}(\overline{M})$ is an $(n+1)$-dimensional
anti-invariant submanifold of $\mathbb{S}^{2n+1}$ invariant under
the flow-action of the characteristic vector field $\xi_0$ and,
locally, $\widetilde{M}$ is isometric to $\mathbb{S}^1\times M^n$.
The submanifold $\overline{M}$ is a parallel Lagrangian submanifold
if and only if $M$ is an integral $\mathcal{C}$-parallel submanifold
(see ~\cite{HN}), and it was proved in ~\cite{DFELSMCO} that a
parallel Lagrangian submanifold $\overline{M}$ is biharmonic if and
only if $M$ is $(-4)$-biharmonic.

We recall here that a map $\psi:(M,g)\to (N,h)$ is
\textit{$(-4)$-biharmonic} if it is a critical point of the
$(-4)$-bienergy $ E_2(\psi)-4E(\psi)$, i.e. $\psi$ verifies
$\tau_2(\psi)+4\tau(\psi)=0$. Also, a real submanifold
$\overline{M}$ of $\mathbb{C}P^n$ is called \textit{Lagrangian} if
it has dimension $n$ and the complex structure $\overline{J}$ of
$\mathbb{C}P^n$ maps the tangent space to $\overline{M}$ onto the
normal one.

Thus, in order to determine all proper-biharmonic parallel
Lagrangian submanifolds of $\mathbb{C}P^3$, we shall determine the
$(-4)$-biharmonic integral $\mathcal{C}$-parallel submanifolds of
$\mathbb{S}^7(1)$.

Just as in the case of Theorem ~\ref{t0} we obtain

\begin{theorem}
The integral submanifold ${\bf i}:M^{3}\to\mathbb{S}^7(1)$ is $(-4)$-biharmonic
if and only if
$$
\begin{cases}\Delta^{\perp}H+\trace
B(\cdot,A_{H}\cdot)
-7H=0\\
4\trace
A_{\nabla^{\perp}_{(\cdot)}H}(\cdot)+3\grad(|H|^{2})=0.\end{cases}
$$

\end{theorem}

Therefore it follows

\begin{proposition} A non-minimal integral $\mathcal{C}$-parallel
submanifold $M^{3}$ of $\mathbb{S}^{7}(1)$ is $(-4)$-biharmonic if
and only if
\begin{equation}\label{eq:-4}
\trace B(\cdot,A_{H}\cdot)=6H.
\end{equation}
\end{proposition}

Now, we can state

\begin{theorem}\label{t3} A $3$-dimensional integral $\mathcal{C}$-parallel submanifold
$M^{3}$ of $\mathbb{S}^{7}(1)$ is $(-4)$-biharmonic if and only if
either
\begin{enumerate}
\item $M^{3}$ is flat and locally is a product of three curves:
\begin{itemize}
\item The $X_{1}$-curve is a helix with curvatures
$\kappa_{1}=\frac{\lambda^{2}-1}{\lambda}$ and $\kappa_{2}=1$,

\item The $X_{2}$-curve is a helix of order $4$ with curvatures
$\kappa_{1}=\sqrt{\lambda^{2}+\alpha^{2}}$,
$\kappa_{2}=\frac{\alpha}{\kappa_{1}}\sqrt{\lambda^{2}+1}$ and
$\kappa_{3}=-\frac{\lambda\sqrt{\lambda^{2}+1}}{\kappa_{1}}$,

\item The $X_{3}$-curve is a helix of order $4$ with curvatures
$\kappa_{1}=\sqrt{\lambda^{2}+\gamma^{2}+\delta^{2}}$,
$\kappa_{2}=\frac{\delta}{\kappa_{1}}\sqrt{\lambda^{2}+\gamma^{2}+1}$
and
$\kappa_{3}=\frac{\kappa_{2}\sqrt{\lambda^{2}+\gamma^{2}}}{\delta}$,
if $\delta\neq 0$, or a circle with curvature
$\kappa_{1}=\sqrt{\lambda^{2}+\gamma^{2}}$, if $\delta=0$,
\end{itemize}
\noindent where $\lambda,\alpha,\gamma,\delta$ are constants given
by
\begin{equation}\label{system0-4}
\begin{cases}
(3\lambda^{2}-1)(3\lambda^{4}-8\lambda^{2}+
1)+\lambda^{4}((\alpha+\gamma)^{2}+\delta^{2})=0\\
(\alpha+\gamma)(5\lambda^{2}+\alpha^{2}+\gamma^{2}-7)+\gamma\delta^{2}=0\\
\delta(5\lambda^{2}+\delta^{2}+3\gamma^{2}+\alpha\gamma-7)=0\\
1+\lambda^{2}+\alpha\gamma-\gamma^{2}=0
\end{cases}
\end{equation}
such that $-1<\lambda< 0$,
$0<\alpha\leq\frac{\lambda^2-1}{\lambda}$, $\alpha\geq\delta\geq 0$,
$\alpha>2\gamma$ and $\lambda^{2}\neq\frac{1}{3}$;

\noindent or

\item $M^{3}$ is locally isometric to a product
$\Gamma\times\bar{M}^{2}$, between a helix with curvatures
$\kappa_{1}=\frac{\sqrt{13}-1}{\sqrt{12-3\sqrt{13}}}$ and
$\kappa_{2}=1$, and a $\mathcal{C}$-parallel surface of
$\mathbb{S}^7(1)$, which is locally isometric to the $2$-dimensional
Euclidean sphere with radius $\sqrt{\frac{3}{7-\sqrt{13}}}$.
\end{enumerate}
\end{theorem}
\begin{proof} It is easy to see that the equation \eqref{eq:-4} is
equivalent to the system
\begin{equation}\label{system-4}
\begin{cases}
(\sum_{i=1}^{3}\lambda_{i})(\sum_{i=1}^{3}\lambda_{i}^{2}-6)+(\alpha+\gamma)(\alpha\lambda_{2}+
\gamma\lambda_{3})\\+(\beta+\delta)(\beta\lambda_{2}+\delta\lambda_{3})=0\\
(\sum_{i=1}^{3}\lambda_{i})(\alpha\lambda_{2}+
\gamma\lambda_{3})+(\alpha+\gamma)(2\lambda_{2}^{2}+\alpha^{2}+3\beta^{2}+\gamma^{2}+
\beta\delta-6)\\+\gamma(\beta+\delta)^{2}=0\\
(\sum_{i=1}^{3}\lambda_{i})(\beta\lambda_{2}+
\delta\lambda_{3})+\beta(\alpha+\gamma)^{2}\\
+(\beta+\delta)(2\lambda_{3}^{2}+\delta^{2}+3\gamma^{2}+\beta^{2}+\alpha\gamma-6)=0.
\end{cases}
\end{equation}

In the same way as for the study of biharmonicity, we shall split
the study of this system, as $M^3$ is given by \textbf{Case I} or
\textbf{Case II} of the classification.

\noindent \textbf{Case I.} The system \eqref{system-4} is equivalent
to the system given by the first three equations of
\eqref{system0-4} and, just like in the proof of Theorem ~\ref{t1},
we conclude.

\noindent \textbf{Case II.}
\begin{enumerate}
\item[(1)] It is easy to verify that this case cannot occur in
this setting.

\item[(2)] The second and the third equation of system
\eqref{system-4} are satisfied and the first equation is equivalent
to $3\lambda^{4}-8\lambda^{2}+1=0$, which solutions are
$\lambda^{2}=\frac{4\pm\sqrt{13}}{3}$. Since $\lambda^{2}<1$ it
follows that $\lambda^{2}=\frac{4-\sqrt{13}}{3}$ and this, together
with the classification of the integral $\mathcal{C}$-submanifolds,
lead to the conclusion.
\end{enumerate}
\end{proof}

Using the explicit equation of the $3$-dimensional integral
$\mathcal{C}$-parallel flat submanifolds in $\mathbb{S}^{7}(1)$ (see
\cite{CBDBTK}) we obtain

\begin{corollary} Any $3$-dimensional flat $(-4)$-biharmonic integral $\mathcal{C}$-parallel
submanifold $M^{3}$ of $\mathbb{S}^{7}(1)$ is given locally by
$$
\begin{array}{ll}
x(u,v,w)=&\frac{\lambda}{\sqrt{\lambda^{2}+1}}\exp(\mathrm{i}(\frac{1}{\lambda}u))\mathcal{E}_{1}
+\frac{1}{\sqrt{(\gamma-\alpha)(2\gamma-\alpha)}}\exp(-\mathrm{i}(\lambda
u-(\gamma-\alpha)v))\mathcal{E}_{2}\\
\\&+\frac{1}{\sqrt{\rho_{1}(\rho_{1}+\rho_{2})}}\exp(-\mathrm{i}(\lambda u+\gamma
v+\rho_{1}w))\mathcal{E}_{3}\\
\\&+\frac{1}{\sqrt{\rho_{2}(\rho_{1}+\rho_{2})}}\exp(-\mathrm{i}
(\lambda u+\gamma v-\rho_{2}w))\mathcal{E}_{4},
\end{array}
$$
where
$\rho_{1,2}=\frac{1}{2}(\sqrt{4\gamma(2\gamma-\alpha)+\delta^{2}}\pm\delta)$,
$-1<\lambda<0$, $0<\alpha\leq\frac{\lambda^2-1}{\lambda}$,
$\alpha\geq\delta\geq 0$, $\alpha>2\gamma$,
$\lambda^{2}\neq\frac{1}{3}$, the tuple
$(\lambda,\alpha,\gamma,\delta)$ being one of the following
$$
\Bigg(-\sqrt{\frac{4-\sqrt{13}}{3}}, \ \sqrt{\frac{7-\sqrt{13}}{6}},
\ -\sqrt{\frac{7-\sqrt{13}}{6}}, \ 0\Bigg),
$$
$$
\Bigg(-\sqrt{\frac{1}{5+2\sqrt{3}}}, \
\sqrt{\frac{45+21\sqrt{3}}{13}}, \ -\sqrt{\frac{6}{21+11\sqrt{3}}},
\ 0\Bigg),
$$
or
$$
\Bigg(-\sqrt{\frac{1}{6+\sqrt{13}}}, \
\sqrt{\frac{523+139\sqrt{13}}{138}}, \
-\sqrt{\frac{79-17\sqrt{13}}{138}}, \
\sqrt{\frac{14+2\sqrt{13}}{3}}\Bigg),
$$
and $\{\mathcal{E}_i\}_{i=1}^4$ is an orthonormal basis of
$\mathbb{C}^4$ with respect to the usual Hermitian inner product.
\end{corollary}

\begin{proof}
In order to solve the system \eqref{system0-4}, we first note that,
since $\alpha>2\gamma$, from the fourth equation it results
$\gamma<0$.

The third equation suggests that we need to split our study in two
cases as $\delta$ is equal to $0$ or not.

\noindent\textbf{Case 1: $\delta=0$.} In this case the third
equation holds whatever the values of $\lambda$, $\alpha$ and
$\gamma$ are, and so does the condition $\alpha\geq\delta$.

If $\alpha=-\gamma$ we easily obtain that the solution of the system
is
$$
\lambda=-\sqrt{\frac{4-\sqrt{13}}{3}},\quad\alpha=\sqrt{\frac{7-\sqrt{13}}{6}},
\quad\gamma=-\sqrt{\frac{7-\sqrt{13}}{6}}.
$$

In the following, we shall look for $\alpha$ of the form
$\alpha=\omega\gamma$, where $\omega\in(-\infty,0)\setminus\{-1\}$.
From the second and the fourth equations of the system we have
$\lambda^2=-\frac{\omega^2+7\omega-6}{(\omega-2)(\omega-3)}$
$\gamma^2=\frac{12}{(\omega-2)(\omega-3)}$ and then
$\alpha^2=\frac{12\omega^2}{(\omega-2)(\omega-3)}$. Replacing in the
first equation, after a straightforward computation, it can be
written as
$$
3\omega^6+16\omega^5-58\omega^4-140\omega^3+531\omega^2-444\omega+108=0,
$$
which is equivalent to
$$
(\omega-2)^2(3\omega^4+28\omega^3+42\omega^2-84\omega+27)=0,
$$
whose solutions are $2$, $-3\pm 2\sqrt{3}$ and $\frac{-5\pm
2\sqrt{13}}{3}$. From these solutions the only one to verify the
supplementary conditions is $\omega=-3-2\sqrt{3}$, for which we have
$$
\lambda=-\sqrt{\frac{1}{5+2\sqrt{3}}},
\quad\alpha=\sqrt{\frac{45+21\sqrt{3}}{13}},
\quad\gamma=-\sqrt{\frac{6}{21+11\sqrt{3}}}.
$$

\noindent\textbf{Case 2: $\delta>0$.} In this case the third
equation of \eqref{system0-4} becomes
$$
5\lambda^2+\delta^2+3\gamma^2+\alpha\gamma-7=0.
$$
Now, again taking $\alpha=\omega\gamma$, this time with
$\omega\in(-\infty,0)$, from the last three equations of the system,
we easily get
$\lambda^2=-\frac{\omega^2+9\omega+2}{(\omega-1)(\omega-2)}$,
$\alpha^2=\frac{12\omega^3}{(\omega-1)^2(\omega-2)}$,
$\gamma^2=\frac{12\omega}{(\omega-1)^2(\omega-2)}$ and
$\delta^2=\frac{12(\omega+1)^2}{(\omega-1)^2}$. Replacing in the
first equation of the system we obtain the solutions $-2\pm\sqrt{3}$
and $-4\pm\sqrt{13}$, from which only $\omega=-4-\sqrt{13}$ verifies
the supplementary conditions. Therefore, we obtain
$$
\lambda=-\sqrt{\frac{1}{6+\sqrt{13}}},
\quad\alpha=\sqrt{\frac{523+139\sqrt{13}}{138}},
$$
$$
\gamma=-\sqrt{\frac{79-17\sqrt{13}}{138}},
\quad\delta=\sqrt{\frac{14+2\sqrt{13}}{3}},
$$
and we conclude.
\end{proof}

\begin{remark}
By a straightforward computations we can check that the images of
the cylinders over the above $x$ are, respectively: the Riemannian
product of a circle of radius $\sqrt{\frac{5-\sqrt{13}}{12}}$ and
three circles, each of radius $\sqrt{\frac{7+\sqrt{13}}{36}}$; the
Riemannian product of two circles each of radius
$\sqrt{\frac{3+\sqrt{3}}{12}}$ and two circles each of radius
$\sqrt{\frac{3-\sqrt{3}}{12}}$; the Riemannian product of a circle
of radius $\sqrt{\frac{5+\sqrt{13}}{12}}$ and three circles, each of
radius $\sqrt{\frac{7-\sqrt{13}}{36}}$.
\end{remark}

\end{document}